\documentclass[pdflatex,sn-mathphys-num,]{sn-jnl}

\usepackage{graphicx}%
\usepackage{multirow}%
\usepackage{amsmath,amssymb,amsfonts}%
\usepackage{amsthm}%
\usepackage{mathrsfs}%
\usepackage[title]{appendix}%
\usepackage{xcolor}%
\usepackage{textcomp}%
\usepackage{manyfoot}%
\usepackage{booktabs}%
\usepackage{algorithm}%
\usepackage{algorithmicx}%
\usepackage{algpseudocode}%
\usepackage{listings}%

\theoremstyle{thmstyleone}%
\newtheorem{theorem}{Theorem}
\newtheorem{corollary}{Corollary}

\theoremstyle{thmstyletwo}%
\newtheorem{example}{Example}%
\newtheorem{remark}{Remark}%

\theoremstyle{thmstylethree}%
\newtheorem{definition}{Definition}%

\usepackage{enumitem,mathrsfs}
\usepackage[normalem]{ulem}
\usepackage{mathtools}
\usepackage{ upgreek }

\makeatletter
\newcommand*\bigcdot{\mathpalette\bigcdot@{.5}}
\newcommand*\bigcdot@[2]{\mathbin{\vcenter{\hbox{\scalebox{#2}{$\m@th#1\bullet$}}}}}
\makeatother

\definecolor{ForestGreen}{RGB}{34,139,34}

\newtheorem*{theorem*}{Theorem}
\newtheorem*{definition*}{Definition}
\newtheorem*{conjecture*}{Conjecture}
\newtheorem*{lemma*}{Lemma}
\newtheorem*{example*}{Example}
\newtheorem*{remark*}{Observacion}

\newcommand{\ZZ}{\ensuremath{{\mathbb Z}}}
\newcommand{\CC}{\ensuremath{{\mathbb C}}}
\newcommand{\CW}{\ensuremath{{\widehat{\mathbb C}}}}
\newcommand{\RR}{\ensuremath{{\mathbb R}}}
\newcommand{\NN}{\ensuremath{{\mathbb N}}}
\newcommand{\HH}{\ensuremath{{\mathbb H}^{2} }}
\font\myfont=cmr10 at 12pt
\newcommand{\e}{{\text{\myfont e}}}

\renewcommand{\Re}[1]{{\mathfrak{Re}\left(#1\right)}}
\renewcommand{\Im}[1]{{\mathfrak{Im}\left(#1\right)}}

\newcommand{\del}[2]{\frac{\partial #1}{\partial #2}}

\usepackage{tikz,etoolbox}
\newcommand{\circled}[1]{ \text{\tikz[baseline=(char.base)]{
\node[shape=circle,draw,inner sep=1pt] (char) {$#1$};}} }
\robustify{\circled}

\usepackage[pagewise]{lineno}

\usepackage[normalem]{ulem}

\begin{document}

\title[Tessellations of rational functions \& Riemann's existence thm]{Tessellations of rational complex functions and 
the Riemann's existence theorem}


\author[1]{\fnm{} \sur{Alvaro Alvarez--Parrilla}}\email{alvaro.uabc@gmail.com}
\equalcont{These authors contributed equally to this work.}

\author*[2]{\fnm{} \sur{Roberto Guti\'errez--Soto}}\email{rgutierrez.s@outlook.com}
\equalcont{These authors contributed equally to this work.}

\author[3]{\fnm{} \sur{Jes\'us Muci\~no--Raymundo}}\email{muciray@matmor.unam.mx}
\equalcont{These authors contributed equally to this work.}

\affil[1]{\orgname{Grupo Alximia SA de CV}, \orgaddress{\city{Ensenada}, 
\state{Baja California}, \country{M\'exico}}}

\affil*[2]{\orgname{Universidad Aut\'onoma del Estado de M\'exico}, \orgaddress{\city{Toluca}, 
\state{Estado de M\'exico}, \country{M\'exico}}}

\affil[3]{\orgdiv{Centro de Ciencias Matem\'aticas}, \orgname{Universidad Nacional Aut\'onoma de M\'exico}, \orgaddress{\street{} \city{Morelia}, 
\state{Michoac\'an}, \country{M\'exico}}}


\abstract{A complex rational function 
$R$, of degree $n \geq 2$, on a compact Riemann surface $M$ 
provided with a cyclic order on its $q$ distinct
critical values,
determines an homogeneous tessellation of the Riemann surface $M$, 
whose $2n$ tiles are topological $q$--gons with alternating
colors.  
The tessellation provides a simple and straightforward 
visual description of the rational function $R$.
Conversely, assume a possibly non homogeneous 
tessellation $\mathscr{T}$ of a compact 
oriented $C^1$ surface $\mathcal{M}$, 
with tiles of alternating colors and a suitable labelling 
in the vertices of its tiles. 
Non homogeneous means that the tiles of 
$\mathscr{T}$ are $\rho$--gons, for different values of $\rho$.
Then there exists a Riemann surface
structure $M$ on $\mathcal{M}$, a complex rational function $R$
and a cyclic order of its critical values,
such that the tessellation of $R$ on $M$ 
topologically coincides with the original $\mathscr{T}$.}




\maketitle

\section{Introduction}

The field of complex rational functions on a compact 
Riemann surface $M$ encodes the algebraical, analytical and 
geometrical, aspects of the complex structure of $M$.  
Our starting point is the following 
classical idea, that we call, the Schwarz--Klein's algorithm:

\smallskip

{\it
\centerline{
A complex rational function 
$R: M \longrightarrow \CW_w$, of degree $n \geq 2$,}

\centerline{
provided with a cyclic order  $\mathcal{L}_\gamma$ of its critical values,
}

\centerline{

determines a tessellation $\mathscr{T}_\gamma(R)$ of
the Riemann surface $M$ }

\centerline{
whose $2n$ tiles are topological $q$--gons 
with alternating colors. }

}

\medskip

\noindent 
Here $q \geq 2$ is the number of distinct critical values of $R$.  
Certainly, 
the  algorithm requires a little bit of information;
a \emph{cyclic order $\mathcal{L}_\gamma$}
of the $q$ critical values of $R$. This is provided by
an oriented Jordan path $\gamma$ in $\CW_w$ 
running through the critical values of $R$. 

\noindent
Hence, 
the input data of the algorithm is a triad 
$(M, R, \mathcal{L}_\gamma)$.
The
output is the homogeneous tessellation $\mathscr{T}_\gamma(R)$, 
it allows for a simple/straightforward 
visual description of the rational function $R$.
A homogeneous tessellation is a collection of
copies of two topological $q$--gons\footnote{
The tiles $T$ and $T^\prime$ arise as the complement of 
$\gamma$ in $\CW_w$: that is 
$\CW_w \backslash \gamma = T \cup T^\prime$.},
$T$ and $T^\prime$, 
that fit together to cover the surface $M$, without overlaps or gaps.
Non homogeneous means that the tiles 
are $\rho$--gons, for different values of $\rho$.

\noindent
As far as we known this idea has its roots in the works of
H.\,A. Schwarz and F.\,Klein.  

In search of an inverse result, let us consider a 
possibly non homogeneous
tessellation $\mathscr{T}$ on 
a compact oriented $C^1$ surface $\mathcal{M}$, 
whose $2n \geq 4$ tiles, are topological $\rho$--gons 
with alternating colors.

\smallskip

{\it
\centerline{ 
How can we recognize that $\mathscr{T}$ is 
equivalent to a tessellation $\mathscr{T}_\gamma(R)$, 
up to }

\centerline{
orientation preserving  homeomorphism, 
as in the Schwarz--Klein's algorithm?
}
}

\smallskip

\noindent
Around 2010, 
W.\,P.\,Thurston conducted a group discussion 
on what he called {\it the shape of rational maps}. 
As a result of the group discussion, 
S. Koch {\it et al.} \cite{Koch-Lei} 
considering tessellations of $\mathbb{S}^2$, found specific conditions, namely 
that the underlying graph is globally and locally balanced,
under which the tessellation
can be recognized topologically as 
a generic branched covering.

In fact, Thurston's notion of consistent labelling 
(a consequence of the global and local balance of the 
tessellation) and our
notion of \emph{consistent $q$--labelling} $\mathcal{L}_c$,
provide a labelling of the vertices of the tessellation,
consistent with the labelling of the vertices of the 
individual tiles as $q$--gons, also
taking into account the orientation according to their colouring and 
another technical condition; 
see Definition \ref{def-etiquetado-consistente} 
for full details.

In order to provide an answer to the previous question, 
let $M=(\mathcal{M}, J)$ be a compact Riemann surface,
$\mathcal{M}$ denotes a 
compact oriented $C^1$
surface of genus ${\tt g}$ and 
$J:T \mathcal{M} \longrightarrow T \mathcal{M}$ is
the complex structure, $J\circ J = -Id$. 
The orientation of $\mathcal{M}$ coincides 
with the one induced by $J$.
The main result is as follows.

\begin{theorem}
\label{teo:principal}
Let $M=(\mathcal{M}, J)$ be a compact 
Riemann surface of genus ${\tt g}$, 
and  $n \geq 2$ be an integer number.
\begin{enumerate}[label=\arabic*),leftmargin=*]
\item 
A rational function $R: M \longrightarrow \CW_w$, 
of degree $n$, provided with a cyclic order $\mathcal{L}_\gamma$ of
its $q$ critical values,
determines a  homogeneous tessellation 
$$
\mathscr{T}_\gamma(R) = 
\underbrace{
T_1 \cup \ldots \cup T_n}_{ \text{blue tiles} }
\cup 
\underbrace{ T^\prime_1 \cup  
\ldots  \cup T^\prime_n}_{\text{gray tiles}} 
\subset M ,
$$

\noindent  
whose tiles are topological $q$--gons 
with alternating colors, and
a consistent $q$--labelling $R^*\mathcal{L}_\gamma$.

\smallskip 

\item
Let $\mathscr{T}$ be a possibly non homogeneous 
tessellation of $\mathcal{M}$, 
whose $2n$ tiles 
are topological $\rho$--gons 
with alternating colors.
Assume in addition that $\mathscr{T}$ is 
provided with a consistent $q$--labelling $\mathcal{L}_c$, 
such that
$2+2{\tt g} \leq \rho  \leq q \leq 2n + 2{\tt g} -2$.
Then $(\mathscr{T}, \mathcal{L}_c)$ determine 
a (non unique) Riemann surface $M$,
a rational function

\centerline{$R: M \longrightarrow \CW_w$,} 

\noindent 
and a Jordan path $\gamma$ 
satisfying that the tessellation
$\mathscr{T}_\gamma(R)$ is $\mathscr{T}$, up to
orientation preserving homeomorphism of $\mathcal{M}$.
\end{enumerate}
\end{theorem}

A cyclic order $\mathcal{L}_\gamma$ determines a consistent $q$--labelling $R^*\mathcal{L}_\gamma$
for the given rational function $R$.
Conversely, there are tessellations $\mathscr{T}$ that
admit consistent $q$--labellings for different values of $q$;
see Example \ref{ex:racional-cocriticos}, 
Figures \ref{fig:tres-mosaicos-etiquetas-h}.c and \ref{fig:tres-mosaicos-etiquetas-q=4-5-5}.

The pair $(\mathscr{T}, \mathcal{L}_c)$ in assertion (2) 
gives rise to an analytical triple  
$(M, R,\mathcal{L}_\gamma)$, which in turn determines 
the pair $(\mathscr{T}_\gamma(R), R^*\mathcal{L}_{\gamma})$, 
which is topologically orientation preserving equivalent 
to the original
pair $(\mathscr{T}, \mathcal{L}_c)$.

\noindent
However, 
by starting with an analytical triple 
$(M,R,\mathcal{L}_\gamma)$, it determines
a pair $(\mathscr{T}_\gamma(R), R^*\mathcal{L}_{\gamma})$, 
which in turn gives rise to 
an analytical triple 
$(\widetilde{R},\widetilde{M},\widetilde{\gamma})$ 
which generically is not  analytically
equivalent to the original $(R,M,\gamma)$.

In order to geometrically perform the gluing,
implicit in assertion (2) of Theorem \ref{teo:principal},
we furnish the sphere with a holomorphic vector field 
$(\CW_w, \partial / \partial w )$,
that is with a translation structure,
see Corollary \ref{cor:funciones-racionales-por-cirujia}.
Moreover, 
by showing that the critical values are
parameters involved in the construction,
one can obtain many different triads
$(M, R, \mathcal{L}_\gamma)$ from the same tessellation
$(\mathscr{T}, \mathcal{L}_c)$.

\smallskip

A very brief list of
technical and historical\footnote{
As usual, 
advances in the theory have been discovered 
independently by diverse authors.}
comments is below.

\noindent 
$\bigcdot$
We attribute assertion (1) of Theorem \ref{teo:principal} to 
H.\,A.\,Schwarz 
(recall his study about uniformization of polygonal regions, 
and the resulting tessellations 
\cite{Schwarz}, a contemporary description is in 
\cite[ch.\,5]{Yoshida})
and 
to F.\,Klein following 
\cite{Klein1}, 
\cite{Chislenko-Tschinkel}.
Hence, we introduce the Schwarz--Klein's algorithm in 
\S \ref{sec:algoritmo-Schwarz-Klein-racional}.
A.\,Speiser \cite{Speiser} and
R.\,Nevanlinna \cite{Nevanlinna2}, considered
tessellations for transcendental functions, 
also see R.\,Peretz \cite{Peretz}.

\noindent
$\bigcdot$ 
The case of $q=3$ critical values 
is essentially a first step in the Theory of 
dessins d'enfants, 
due to G.\,V.\,Bely{\u\i} \cite{Belyi} and 
A.\,Grothendieck \cite{Grothendieck}, currently 
a highly developed subject.
See the books 
S.\,K.\,Lando {\it et al.} \cite{Lando-Zvonkin}, 
A.\,Degtyarev \cite{Degtyarev} and
G.\,A.\,Jones {\it et al.} \cite{Jones-Wolfart} remarkable for its historic review.

\noindent 
$\bigcdot$ 
Assertion (1) of Theorem \ref{teo:principal} is 
a statement, 
using tessellations, of the Riemann's existence theorem, compare with 
the versions in
I.\,Bauer {\it et al.} \cite{Bauer-Catanese},
H. Volklein \cite[Ch.\,6]{Volklein},
S.\,I.\,Lando  {\it et al.} \cite[p.\,74]{Lando-Zvonkin}
and
R.\,Cavalieri {\it et al.} \cite[p.\,85]{Cavalieri-Miles}.

\noindent 
$\bigcdot$ 
The role of the Jordan path $\gamma$
in assertion (1) is to provide a cyclic order $\mathcal{L}_\gamma$
to the critical values of $R$;  see 
Corollary  \ref{cor:el-orden-ciclico-para-los-valores}.
On $M=\CW_z$, the change of the order was studied by
H.\,Habsch \cite{Habsch}, 
see also S.\,Lando {\it et al.} \cite[5.4]{Lando-Zvonkin}
for the polynomial case.

\noindent 
$\bigcdot$
The proof of assertion (2) can be described succinctly
as ``pullback the complex structure'', 
compare with S.\,Stoilow \cite{Stoilow}.
However, this road requires heavy analytic machinery.
Our proof uses elementary surgery techniques, as
Corollary \ref{cor:funciones-racionales-por-cirujia}
emphasizes.
In fact, it can be widely applied to 
complex analytic: 
functions,
quadratic differentials 
K.\,Strebel \cite[\S\,12.3]{Strebel} 
(we learned this surgery technique from him), 
H.\,G.\,Dias--Marin \cite{Diaz}, and
vector fields 
J.\,Muci\~no--Raymundo \cite{Mucino}.

\noindent
$\bigcdot$
In the case of $M=\CW_z$, the rational functions sharing 
a set of critical values 
was described algebraically by J.\,Mycielski \cite{Mycielski}, 
compare with the classical result of L.\,R.\,Goldberg
\cite{Goldberg}.
Recently, K.\,Lazebnik \cite{Lazebnik} provides
an analytic study of the paths $R^{-1}(\gamma)$.
 
\noindent
$\bigcdot$
Rational functions with $q$--critical values
give rise to subfamilies of $F(w_1, \ldots , w_q)$, the
Riemann surfaces
whose branch points lie over a finite set of singular values 
$\{ w_1, \ldots, w_q\}  \subset \CW_w$, where finite and/or 
infinite branchings are
allowed; see R.\,Nevanlinna \cite[ch.\,XI,\S2]{Nevanlinna2}.
In order to describe a meromorphic function $w(z)$, whose
Riemann surface belongs to $F(w_1, \ldots , w_q)$, 
the usual tool is its {\it Speiser graph}. 
See for instance 
A.\,A.\,Goldberg {\it et al.} \cite[ch.\,7\,\S4]{GoldbergOstrovskii}.
In the topological sense, a Speiser graph is the dual of a tessellation 
associated to $w(z)$.  

\noindent
$\bigcdot$
Tessellations have been applied to the study 
of differential equations, see 
D.\,Masoero \cite{Masoero}, 
G.\,Le\'on--Gil {\it et al.} \cite{Leon-Mucino},
A.\,Alvarez--Parrilla {\it et al.} \cite{AlvarezMucino3}, 
and references therein.

\smallskip

The article is organized as follows. 
In Section \ref{sec:Tessellations-and-graphs-on-Riemann-surfaces}, 
we recall the definition of tessellation
and its associated graphs. The Schwarz--Klein's algorithm
is introduced in \S \ref{sec:algoritmo-Schwarz-Klein-racional}.
Section \ref{sec:Converse-Schwarz-Klein-algorithm} contains the proof 
of assertion (2) of Theorem \ref{teo:principal}.
In \S\,3.1, we review the relationship between analytic and combinatorial
structures involved in the work.
A construction of rational functions 
by surgery is provided in Corollary 
\ref{cor:funciones-racionales-por-cirujia}. 
Section \ref{sec:Examples} presents examples and a brief comment on future directions.

\section{Tessellations and graphs on Riemann surfaces}
\label{sec:Tessellations-and-graphs-on-Riemann-surfaces}

Tessellations and graphs  
appear in many instances in the study of
complex analytic functions and Riemann surfaces,  
with very intricate meanings and notations.
We provide accurate ad hoc concepts. 

\noindent 
Let $M=(\mathcal{M}, J)$ be a compact Riemann surface,
as in the Introduction.
Definitions 
\ref{definicion-de-teselacion}--\ref{def-etiquetado-consistente}
below are considered over the $C^1$ surface 
$\mathcal{M}$, hence they apply on $M$.

\begin{definition}
\label{definicion-de-teselacion}
\begin{upshape}
A {\it tessellation} of a
surface $\mathcal{M}$ is a collection 
with alternating colors
\begin{equation}
\label{teselacion}
\mathscr{T}
=
\underbrace{
T_1 \cup \ldots \cup T_n}_{ \text{blue tiles} }
\cup 
\underbrace{ T^\prime_1 \cup  
\ldots  \cup T^\prime_n}_{\text{gray tiles}} \subset 
\mathcal{M}, 
\ \ \ 
n \geq 2,
\end{equation}

\noindent
where the $2n$ {\it tiles}
$\{ T_{\alpha}, \, T_{\alpha}^\prime \}_{\alpha=1}^n$
are open Jordan domains,
such that: 

\noindent 
i)
The union of their closures 
$\cup_{{\alpha}=1}^{n}
\big( \overline{T_{\alpha}} 
\cup
\overline{T_{\alpha}^\prime} \big)$
is $\mathcal{M}$.

\noindent 
ii)
If the intersection of the closures of any
two tiles is non--empty, 
then it consists of
a finite number of 
simple paths (edges) and their extreme points (vertices).
\end{upshape}
\end{definition}

A  tessellation $\mathscr{T}$
has $n$ blue tiles and $n$ gray tiles, this is 
called the {\it global balance condition} in
\cite{Koch-Lei}.  
By looking at the boundaries of the tiles,
say $\partial \overline{T_\alpha}, \,  
\partial \overline{T_{\alpha}^\prime}$,
a tessellation $\mathscr{T}$ 
determines an underlying graph $\Gamma$.

\begin{definition}
\label{definicion-de-t-graph}
\begin{upshape}
A {\it ${\tt t}$--graph} $\Gamma$ is a finite oriented 
connected graph 
embedded in $\mathcal{M}$, 
with vertices $V(\Gamma)$ of even valence 
equal or greater than $4$  
and edges $E (\Gamma)$, such that:

\noindent 
i) \ $\mathscr{T}(\Gamma) \doteq
\mathcal{M} \setminus \Gamma $
is a tessellation,   
as in Definition \ref{definicion-de-teselacion}. 

\noindent 
ii) Each blue tile $T_\alpha$ is on
the left side of the oriented edges of $\Gamma$.
\end{upshape}
\end{definition}

With the above in mind, 
a possibly non homogeneous tessellation $\mathscr{T}$ and a 
${\tt t}$--graph $\Gamma$ are essentially 
equivalent objects, where the alternating 
colouring in Equation 
\eqref{definicion-de-teselacion} corresponds 
to the orientation of the edges in Definition 
\ref{definicion-de-t-graph}. 
In simple words, 
a ${\tt t}$--graph must be understood as the simplest
oriented graph describing a tessellation. 
The tessellations arising from complex rational 
functions require a more accurate notion, as follows.

\begin{definition}
\label{definicion-de-R-map}
\begin{upshape}
An {\it ${\tt R}$--map} $\widehat{\Gamma}$ is a 
finite, oriented, connected graph 
embedded in  $\mathcal{M}$, 
with vertices $V(\widehat{\Gamma})$ of
even valence  equal or greater than $2$
and edges $E (\widehat{\Gamma})$, such that:

\noindent 
i) The subset of vertices of valence equal or greater
than 4 is non empty. 

\noindent 
ii) If we forget the vertices
of valence 2 of $\widehat{\Gamma}$, then
we obtain a ${\tt t}$--graph 
$\Gamma$ such that:

\centerline{
$\mathscr{T}(\widehat{\Gamma}) \doteq 
\mathscr{T}(\Gamma)=
T_1 \cup \ldots \cup T_n \cup 
T^\prime_1 \cup \ldots  \cup T^\prime_n $
}

\noindent
is, set theoretically, a tessellation as in Definition 
\ref{definicion-de-teselacion}.

\noindent 
iii) 
The boundary 
$\partial \overline{T_\alpha}$ 
(resp. $\partial \overline{T_{\alpha}^\prime}$) of  
each tile has $q \geq 2$ edges of $\widehat{\Gamma}$,
{\it i.e.} the tessellation of $\widehat{\Gamma}$ is homogeneous.
\end{upshape}
\end{definition}

The {\it forgetting vertices operation}
in part (ii) above is as follows. 
We consider a vertex $z_1=0$ of valence 2 
in $\widehat{\Gamma}$ and its two adjacent edges, thus we have
$(-1,0) \cup \{0\} \cup (0,1)$.
The operation of forgetting
the vertex $0$ replaces the above by an unique
edge $(-1,1)$.

After Theorem \ref{teo:principal}, the
name ${\tt R}$--map for $\widehat{\Gamma}$ 
must be understood as a 
coarse abbreviation of
``complex rational function''.

\begin{remark}
\label{grado-y-gonalidad}
\begin{upshape}
An ${\tt R}$--map $\widehat{\Gamma}$ 
has two numerical attributes $n, \, q \geq 2$:

\noindent 
$\bigcdot$ 
$n$ is the {\it degree of $\widehat{\Gamma}$}   
and

\noindent 
$\bigcdot$
$q$ is the {\it $q$--gonality}
(the tiles of the tessellation $\mathscr{T} (\widehat{\Gamma})$ 
are topological $q$--gons). 
\end{upshape}
\end{remark}

\begin{example}
Figure \ref{fig:tres-mosaicos-etiquetas-h}.a 
illustrates an affine view of a ${\tt t}$--graph $\Gamma \subset \CW_z$,
whereas 
Figure \ref{fig:tres-mosaicos-etiquetas-h}.c sketches
an ${\tt R}$--map $\widehat{\Gamma} \subset \CW_z$, 
where $\infty \in \CW_z$ is a vertex of both graphs. 
In both cases, the tessellation 
is set theoretically the same;
the difference is on the vertices of valence 2.
As a consequence, 
the tessellation $\mathscr{T}(\Gamma)$ is non 
homogeneous, while $\mathscr{T}(\widehat{\Gamma})$ 
is homogeneous.
Several tessellations $\mathscr{T}(\Gamma)$
are illustrated in \cite{ALRJC} and \cite{Gonzalez-Mucino}.
\end{example}

A main combinatorial feature of our graphs is the following. 

\begin{definition}
\label{def-etiquetado-consistente}
A {\it consistent $q$--labelling 

\centerline{$\mathcal{L}_c: V(\Gamma) \longrightarrow \ZZ_q, 
\ \ \
q \geq 2$, } 

\noindent 
for a ${\tt t}$--graph $\Gamma$} satisfies the following conditions:

\noindent
i) For each blue tile $T_\alpha$ of the tessellation
$\mathscr{T}(\Gamma)$, 
if $\{ z_\iota \}$ are the vertices of its
boundary $\partial \overline{T_\alpha}$, ordered 
with cyclic anti--clockwise sense, then the 
labels $\{ \mathcal{L}_c ( z_\iota ) \} \subset \ZZ_q$
appear exactly once and with the same cyclic 
order as in $\ZZ_q$.

\noindent 
ii) Each label $j \in \ZZ_q$ appears under
$\mathcal{L}_c$ for at least one
vertex $z_\iota \in V(\Gamma)$ 
of valence equal or greater than $4$.  
\end{definition}

\begin{remark}[Consistent $q$--labelling for ${\tt R}$--maps]
\noindent 1. 
The notion of consistent $q$--labelling for 
a ${\tt t}$--graph $\Gamma$ 
extends to any ${\tt R}$--map
$\widehat{\Gamma}$ as follows:
for $\widehat{\Gamma}$ all the labels 
of $\ZZ_q$ appear on the vertices   
of each blue tile $T_\alpha$, since all the tiles of 
$\mathscr{T}(\widehat{\Gamma})$ are $q$--gons.

\noindent 2.
On the other hand, for $\Gamma$ usually some labels
of $\ZZ_q$ are hidden  in the boundary
of each blue tile, since the tiles of the tessellation 
$\mathscr{T}(\Gamma)$ can be 
$\rho$--polygons, for $2 \leq \rho \leq q $,
as we will show
in Example \ref{ex:racional-cocriticos} and Figure 
\ref{fig:tres-mosaicos-etiquetas-h}.
By abuse of notation, we use the notion of consistent 
$q$--labelling for ${\tt t}$--graphs and ${\tt R}$--maps. 

\noindent 3.
Note that not all ${\tt R}$--maps have consistent $q$--labellings.
See for instance figure 10 of \cite{Koch-Lei}.
The figure depicts an ${\tt R}$--map with degree 4
(number of tiles of each color), 
6--gonality (each tile is a 6--gon),
and a labelling in $\ZZ_6$.  
However, all the vertices labelled 5 have valence 2, thus,

\noindent
$\bigcdot$ 
the labelling is not a consistent $6$--labelling,

\noindent
$\bigcdot$ 
the vertices labelled 5 are fake 
cocritical points, equivalently, 
5 is a fake critical value.

\noindent
By forgetting the vertices labelled 5, we obtain another ${\tt R}$--map with degree 4 and
5--gonality, renaming label 6 to label 5 provides a consistent 
$5$--labelling. 

\end{remark}

\subsection{Schwarz--Klein's algorithm for the construction of tessellations}
\label{sec:algoritmo-Schwarz-Klein-racional}

We follow the classical works of 
H.\,A.\,Schwarz \cite{Schwarz}, 
F.\,Klein \cite{Klein1}, \cite{Chislenko-Tschinkel},
R.\,Nevanlinna \cite{Nevanlinna2}\,ch.\,XI\,\S2. 

\smallskip

{\bf Schwarz--Klein's algorithm.}

\noindent 
Let $R:M \longrightarrow \CW_w$ be a rational function 
of degree $n\geq 2$.
Recall the Riemann--Hurwitz formula 

\centerline{$2{\tt g} -2 = -2n + \sum_{z_\iota }
( \mu_\iota -1 )$, }

\noindent 
where $\{z_\iota\}$ are the singular points of $R$ 
with ramification orders $\mu_\iota \geq 2$. 
In particular, if $m$ is the number of critical points of $R$, 
then \, $2+2{\tt g} \leq m \leq  2n +2{\tt g}-2 $.

\noindent 
{\it Step 1.}
Compute the critical points of $R$, 

\centerline{$ 
\mathcal{CP}_R =
\{z_1, \ldots, z_m \} \subset M
$}

\noindent 
and the critical values of $R$, 

\centerline{$
\mathcal{CV}_R =
\{w_1, \ldots, w_j , \ldots , w_q \} \subset \CW_w,
\ \ \ 
j\in 1, \ldots , q, \ 
2+ 2{\tt g} \leq q \leq m$.
}

\noindent 
{\it Step 2.}
Construct an oriented Jordan path 

\centerline{ 
$\gamma \subset \CW_w$ \  running through 
\ $\{w_1, \ldots, w_q \}$,
}

\noindent 
here the cyclic subindices $j\in \ZZ_q$ 
of the critical values are provided by $\gamma$.
We have a trivial tessellation 
\begin{equation}
\label{ec:teselacion-trivial}
\mathscr{T}( \gamma ) = 
\CW_w \backslash \gamma = T \cup T^\prime
\end{equation} 

\noindent
with topological $q$--gons as tiles, 
$T$ is blue (in the left side of $\gamma$) and 
$T^\prime$ is gray.
Moreover, $\gamma$ is a cyclic graph with  $q$ vertices
$\{w_j\} = \mathcal{CV}_R$ and  $q$ edges 
(the respective segments 
$\overline{w_\iota w_{\iota+1}}$, 
$\overline{w _q w_1}$ of $\gamma$).

\noindent 
{\it Step 3.} 
Compute the inverse image of $\gamma$,
\begin{equation}
\label{ec:Gamma}
\Gamma = R^{-1}(\gamma) \subset M.
\end{equation}

\noindent
More accurately, the pullback graph 
\begin{equation}
\label{ec:Gamma-gorrito}
\widehat{\Gamma}=R^* \gamma
\end{equation}

\noindent 
is well defined. 
By definition, a 
{\it cocritial point $\zeta_\kappa \in \mathcal{C}c_R$ 
of $R$} 
satisfies that $R(\zeta_\kappa)$ is a critical 
value but $\upzeta_\kappa$ is a regular point of $R$.
The vertices of $\widehat{\Gamma}$ are 
\begin{equation}
\label{eq-puntos-criticos y-cocriticos}
V(\widehat{\Gamma})
=
\{ \underbrace{\hbox{critical points } \mathcal{CP}_R}_{
\text{even valence } \geq 4 } \} 
\cup 
\{ \underbrace{ \hbox{cocritical points } \mathcal{C}c_R}_{\text{valence } = 2} \} .
\end{equation}

\noindent 
The edges $E(\widehat{\Gamma})$ are 
the respective segments in the ${\tt R}$--map  
$\widehat{\Gamma} = R^* \gamma$.
Set theoretically $\Gamma = \widehat{\Gamma}$,  
however they are isomorphic graphs if and only 
the cocritical point set $\mathcal{C}c_R$ 
is empty.

\noindent 
{\it Step 4.} 
The tessellation determined by $R$ and $\gamma$ is 

\centerline{
$
\mathscr{T}_\gamma(R) =
M \backslash \widehat{\Gamma} =
\underbrace{
T_1 \cup \ldots \cup T_n}_{ \text{blue tiles} }
\cup 
\underbrace{ T^\prime_1 \cup  
\ldots  \cup T^\prime_n}_{\text{gray tiles}}, 
\ \ \ 
n \geq 2
$.
}

\noindent 
{\it Step 5.} 
The cocritical points $\mathcal{C}c_R$ plays a 
crucial role, hence the tiles of
$\mathscr{T}(\widehat{\Gamma})$
are topological $q$--gons.
Finally, we add labels to the 
vertices of $\gamma$ with the labelling map 
\begin{equation}
\label{ec:etiquetado-de-gamma}
\mathcal{L}_\gamma: V(\gamma ) 
\longrightarrow \ZZ_q,
\ \ \
w_j  \longmapsto  j .
\end{equation}

\noindent 
The pullback of the labels determines
a \emph{consistent $q$--labelling for $\widehat{\Gamma}$}
\begin{equation}
\label{ec:etiquetado-de-Gamma}
\begin{array}{rcl}
R^*\mathcal{L}_\gamma: 
\mathcal{CP}_R \cup \mathcal{C}c_R 
&\longrightarrow & \ZZ_q
\\
z_\iota & \longmapsto & \mathcal{L}_\gamma(R(z_\iota)) \, .
\end{array}
\end{equation}

\noindent 
Summing up,  
the input of the Schwarz--Klein's algorithm is the pair 
$(M, R, \mathcal{L}_\gamma)$.
The (equivalent) outputs are:

\noindent
$\bigcdot$
The homogeneous tessellation 
$\mathscr{T}_\gamma(R)$ of $M$ with $2n$ tiles, 
which are topological $q$--gons having alternating
colors.

\noindent 
$\bigcdot$
The ${\tt R}$--map
$\widehat{\Gamma}= R^* \gamma \subset M$
with a consistent $q$--labelling $ R^* \mathcal{L}_\gamma $.

\begin{example}[Tessellation of a non generic rational function]
\label{ex:racional-cocriticos}
\upshape
Let 

\centerline{$
R(z):\CW_z \longrightarrow \CW_w,   
\ \ \  
z \longmapsto {z(z^2-1)(z^2-4)}/{(z-3)} 
$}

\noindent
be a non generic\footnote{A rational function 
$R(z)$ in $\CW_z$ of degree 
$n$ is {\it generic} when it has 
$2n-2$ distinct critical points. 
} 
rational function of degree $5$,
as a straightforward computation shows, 
there are six real critical points 
$\{z_1, \ldots z_5, \infty \}$, 
with five real critical values and infinity.
Hence the choice of 
$\gamma= \RR \cup \{ \infty \}$ is suitable.
Figure \ref{fig:tres-mosaicos-etiquetas-h}.a sketches
the affine view of $\Gamma$, the finite critical points 
and its tessellation $\mathscr{T} (\Gamma)$.
The path $\gamma$ induces the cyclic
order $w_1 < w_2 <w_3 < w_4 < w_5 < w_6 =\infty$.
By Equation \eqref{ec:etiquetado-de-gamma},
there is a labelling for the critical values of $R(z)$, 

\centerline{$
\mathcal{L}_\gamma: \mathcal{CV}_R \longrightarrow 
\ZZ_6,
\ \ \
w_j \longmapsto j,
$}

\noindent 
where by simplicity the subindex $j$ coincides with a
cyclic order in $\mathcal{CV}_R$.
Without loss of generality, we define that 
$\mathcal{L}_\gamma(\infty)=6$.

\noindent 
The pullback of $\mathcal{L}_\gamma$ under $R$
determines a consistent $6$--labelling of the 
critical and cocritical points of $R$, that is

\centerline{$
\mathcal{L}_c \doteq 
R^*\mathcal{L}_\gamma: \mathcal{CP}_R \cup \mathcal{C}c_R 
\longrightarrow \ZZ_6, 
\ \ \
z_\iota \longmapsto \mathcal{L}_\gamma(R(z_\iota)).
$}

\noindent 
Figure \ref{fig:tres-mosaicos-etiquetas-h}.c sketches
the affine view of 
$\widehat{\Gamma}= 
R^* \gamma$, 
its homogeneous tessellation $\mathscr{T} (\widehat{\Gamma})$ and the  
consistent $6$--labelling 
$\mathcal{L}_c: 
V(\widehat{\Gamma}) \longrightarrow \ZZ_6$, with 
$\mathcal{L}_\gamma (\infty) =6$. 
In particular, the labelling $\mathcal{L}_c$ 
of $\Gamma$ is computed as

\centerline{$
\begin{array}{rcr}
z_1 \longmapsto R(z_1)=w_5 \longmapsto \mathcal{L}_\gamma(w_5)=5, 
&  \ \ \ \ &
z_4 \longmapsto R(z_4)=w_2 \longmapsto \mathcal{L}_\gamma(w_2)=2, 
\\
z_2 \longmapsto R(z_2)=w_3 \longmapsto  \mathcal{L}_\gamma(w_3) =3, 
&& 
z_5 \longmapsto R(z_4)=w_1 \longmapsto \mathcal{L}_\gamma(w_1)=1,
\\
z_3 \longmapsto R(z_4)=w_4 \longmapsto \mathcal{L}_\gamma(w_4)=4, 
&&
\
z_6=\infty \longmapsto R(z_6)=\infty \longmapsto 
\mathcal{L}_\gamma(w_6)=6.
\end{array}
$}

\noindent 
The $q$--gonality of $\widehat{\Gamma}=R^* \gamma$ is 6, 
{\it i.e.} the tiles are topological $6$--gons. 

\noindent 
The output of the Schwarz--Klein's algorithm applied to $R$ is 
$\mathscr{T}_\gamma (R)$, 
or 
$\widehat{\Gamma}$ with  a labelling $\mathcal{L}_c$
with values in $\ZZ_6$.

\begin{figure}[h!]
\begin{center}
\includegraphics[width=\textwidth]{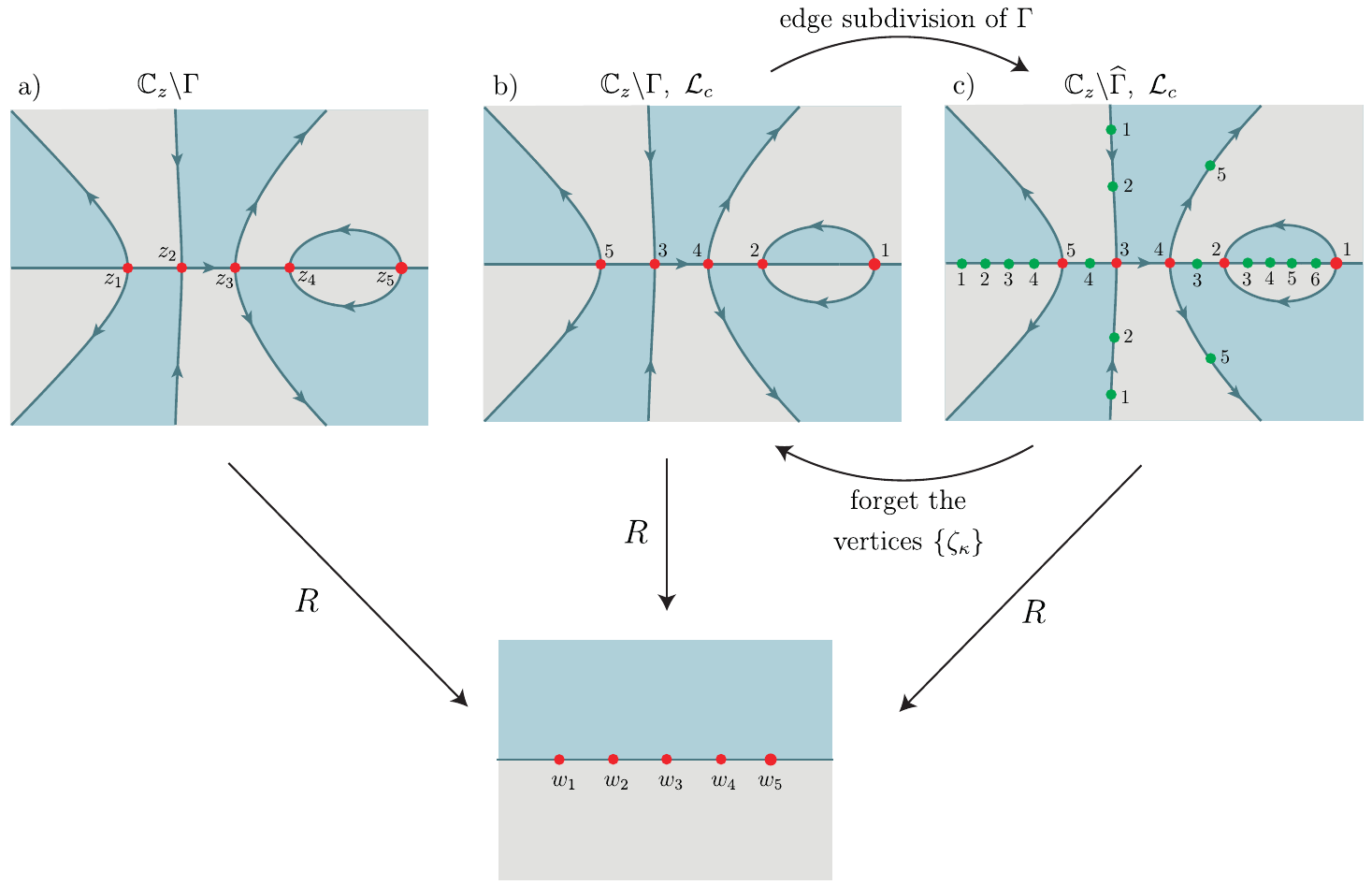}
\caption{
Affine view of the tessellation of a 
non generic rational function $R(z)$ of degree 5. 
a) 
The ${\tt t}$--graph $\Gamma= R^{-1} (\RR \cup \{ \infty\})$ and 
its non homogeneous tessellation
$\mathscr{T}(\Gamma)$. 
b) 
The consistent $6$--labelling 
$\mathcal{L}_c: V(\Gamma) 
\longrightarrow \ZZ_6$, where $\mathcal{L}(\infty)= 6$. 
c) 
The ${\tt R}$--map $\widehat{\Gamma}= R^* \gamma$ and
its homogeneous tessellation
$\mathscr{T}(\widehat{\Gamma})$:
it has a
consistent $6$--labelling 
$R^*\mathcal{L}_c$, each tile 
is a $6$--gon, with
vertices at the (red) critical points $\mathcal{CP}_R$, 
the point $w_6=\infty$ (which has label 6), and
the (green) cocritical points $\mathcal{C}c_R$.
}
\label{fig:tres-mosaicos-etiquetas-h}
\end{center}
\end{figure}

\noindent
As a second part of the example, assume that our departure point
is the ${\tt t}$--graph $\Gamma$ describing the 
non homogeneous 
tessellation $\mathscr{T}(\Gamma)$ 
in Figure \ref{fig:tres-mosaicos-etiquetas-h}.a, 
with $10=2n$ tiles. 
We want to ask about the existence of a rational function $R$
of degree $n=5$,
which topologically realizes the tessellation of $\Gamma$. 
\noindent

\noindent 
As a main difficulty,
note that the tiles of $\mathscr{T}(\Gamma)$ are 
$\rho$--gons, for $\rho=2, \, 3, \, 4$, providing a lower
bound for $q$, on the other hand,
the Riemann--Hurwitz formula provides an upper bound; 
we must therefore search for 
a consistent $q$--labelling with 

\centerline{ 
$\max \{ \rho \} = 4 \leq  q \leq  8 
= 2n  + 2{\tt g} -2 $.
}

\noindent
Further, since ${\tt g}=0$ and noting that condition (ii) of Definition \ref{def-etiquetado-consistente} 
implies $q\leq \# \{$vertices of valence $\geq 4 \} = 6$, or search reduces to $4\leq q\leq 6$.

\noindent 
The key assumption is whether we know 
of consistent $q$--labellings
for the vertices of $\Gamma$. 
In Figure \ref{fig:tres-mosaicos-etiquetas-h}.b, a consistent 6--labelling is shown for the ${\tt t}$--graph.
In order to construct an ${\tt R}$--map that encodes a topological branched covering $\mathcal{R}$ 
over $\CW_w$, we proceed by an edge subdivision operation in $\Gamma$, as in Figure \ref{fig:tres-mosaicos-etiquetas-h}.c.
Finally, 
the existence of a complex structure, making it a rational map
is assertion (2) of Theorem \ref{teo:principal}.

\noindent
As it turns out, for this particular non homogeneous tessellation $\mathscr{T}(\Gamma)$
there are also consistent $q$--labellings for $q=4$, and $q=5$:
Figure \ref{fig:tres-mosaicos-etiquetas-q=4-5-5} shows one for $q=4$ and two for $q=5$;
each gives rise to a different rational function $R$.

\begin{figure}[h!]
\begin{center}
\includegraphics[width=\textwidth]{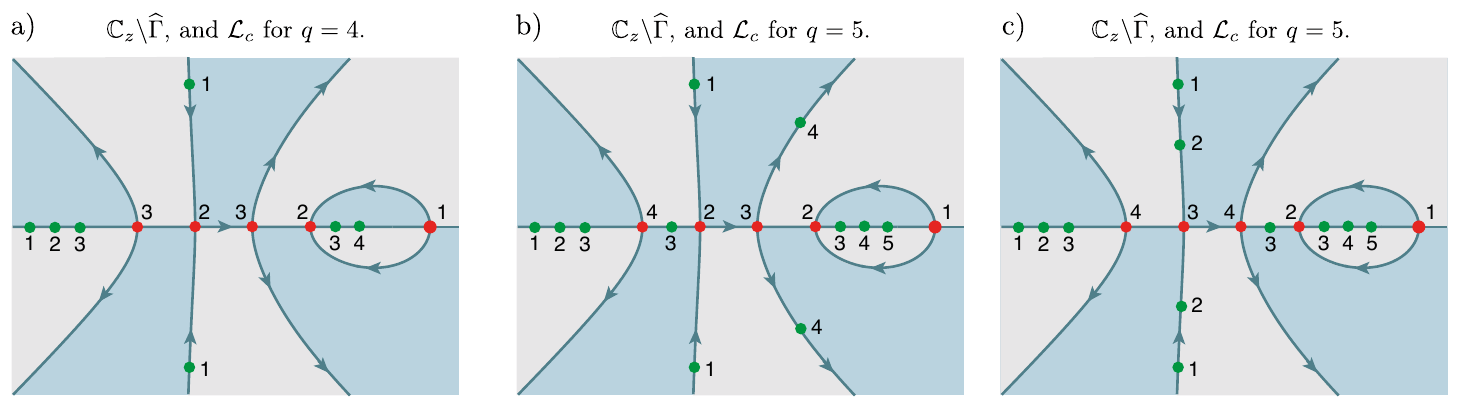}
\caption{
Consistent $q$--labellings for $q=4$, and $5$ 
of the same ${\tt t}$--graph of 
Figure \ref{fig:tres-mosaicos-etiquetas-h}.a
which originates from $R(z)$ in Example 
\ref{ex:racional-cocriticos}.\,
a) 
A consistent $4$--labelling 
$\mathcal{L}_c: V(\Gamma) 
\longrightarrow \ZZ_4$, where $\mathcal{L}(\infty)= 4$.
b)
A consistent $5$--labelling 
$\mathcal{L}_c: V(\Gamma) 
\longrightarrow \ZZ_5$, where $\mathcal{L}(\infty)= 5$. 
c) 
A topologically different 
consistent $5$--labelling 
$\mathcal{L}_c: V(\Gamma) 
\longrightarrow \ZZ_5$, where $\mathcal{L}(\infty)= 5$.
Each case produces a different
${\tt R}$--map $\widehat{\Gamma}= R^* \gamma$ with a
homogeneous tessellation
$\mathscr{T}(\widehat{\Gamma})$.
Each tile is a $q$--gon, with
vertices at the (red) critical points $\mathcal{CP}_R$, 
the point $\infty$ with label $q$, and
the (green) cocritical points $\mathcal{C}c_R$.
}
\label{fig:tres-mosaicos-etiquetas-q=4-5-5}
\end{center}
\end{figure} 
\end{example}

Figure 3.a--b  in  \cite{Gonzalez-Mucino} illustrates
the tessellation $\mathscr{T}(\Gamma) \subset \mathbb{S}^2$ 
of a ${\tt t}$--graph that admits consistent
$q$--labellings for $q=3$ and $4$.

\section{A converse for the Schwarz--Klein's algorithm}
\label{sec:Converse-Schwarz-Klein-algorithm}
Let $\Gamma \subset \mathcal{M}$ be 
a ${\tt t}$--graph, 
where $\mathcal{M}$ is a compact oriented $C^1$ surface  
and
such that $\mathscr{T}(\Gamma) =
\mathcal{M} \backslash \Gamma$ is a tessellation,
as in Definition \ref{definicion-de-teselacion}.

\smallskip 

{\bf Proof of Theorem \ref{teo:principal} assertion (1). }
Follows directly
from the Schwarz--Klein's  algorithm.
Let us comment that this assertion assumes
that $M$ is a Riemann surface, whereas 
assertion (2) will work on $\mathcal{M}$:

\begin{remark}
\label{rem:q=2-3-4}
For $q=2$ critical values, the only functions that appear
are rational functions
of degree $n=2$  on $\CW_z$ with two critical points.

\noindent 
For $q=3$ critical values the pair 
$(\Gamma, \mathcal{L}_c)$ determines a unique complex structure on
$\mathcal{M}$, see \cite[p.\,4]{Jones-Wolfart} 
for an historical comment. 
This fortunate property is key for 
the theory of Bely{\u \i}'s functions and dessins d'enfants. 

\noindent
For  $q\geq 4$, clearly,  a pair $(\Gamma, \mathcal{L}_c)$
determines families of non biholomorphic 
complex structures $(M,J)$ and 
rational functions $R$, see Example \ref{ex:funcion-wp}.1. 
\end{remark}

\smallskip

{\bf Proof of Theorem \ref{teo:principal} assertion (2).}
We perform several steps.

{\it Step 1.}
By using the consistent $q$--labelling 
$\mathcal{L}_c$ and edge subdivision operation 
for $\Gamma$, we get an associated ${\tt R}$--map 
$\widehat{\Gamma}$ in $\mathcal{M}$ with a 
consistent $q$--labelling 
$\mathcal{L}_c: V(\widehat{\Gamma}) \longrightarrow \ZZ_q$, 
as follows.

\noindent
{\it Edge subdivision operation.}
Let $\overline{z_\iota z_\sigma}$ an edge of $\Gamma$ with 
labels 
$\mathcal{L}_c(z_\iota)$ and 
$\mathcal{L}_c(z_\sigma)$. 

\noindent 
If $\mathcal{L}_c(z_\sigma) -\mathcal{L}_c(z_\iota) = 
1  \pmod{q} $, 
then $\overline{z_\iota z_\sigma}$ is an edge of 
$\widehat{\Gamma}$.

\noindent 
If $\mathcal{L}_c(z_\sigma) -\mathcal{L}_c(z_\iota) = 
\nu +1 \geq 2  \pmod{q} $, 
then we consider $\nu$ new vertices in the original edge 
$\overline{z_\iota z_\sigma}$, obtaining $\nu +1$ edges 

\centerline{$ 
\overline{z_\iota \zeta_1}, \ 
\overline{\zeta_1 \zeta_2}, \ldots ,
\overline{\zeta_\nu z_\sigma}
$}

\noindent 
of $\widehat{\Gamma}$. Moreover, the labels of these new 
vertices of valence 2 of $\widehat{\Gamma}$ are

\centerline{$ 
\mathcal{L}_c(z_\iota) = h , \
\mathcal{L}_c(\zeta_1) = h+1 , \
\ldots ,
\mathcal{L}_c(\zeta_1) = h+ \nu , \
\mathcal{L}_c(z_\sigma) = h+\nu+1 \ ;
\ \ 
\pmod{q}.
$}

\noindent 
Figure \ref{fig:tres-mosaicos-etiquetas-h}.b--c illustrates the 
edge subdivision operation.

{\it Step 2.}
We consider the circle 
$\gamma = \RR\cup \{ \infty \} \subset \CW_w$ 
furnished with vertices 

\centerline{$
V(\gamma)=\{ 1, 2,  \ldots, {k-1}, \infty \} $}

\noindent 
and the respective  $q$ segments as edges $E(\gamma)$. 
Thus $\gamma$ determines a trivial tessellation
\eqref{ec:teselacion-trivial}, say

\centerline{
$
\CW_w \backslash \gamma = \HH_+ \cup \HH_- \, , 
$}

\noindent
where as usual $\HH_+$ denotes the open upper half plane.
The labelling of the vertices of $\gamma$ is
$$
\mathcal{L}_\gamma : \{
\underbrace{1,2, \ldots, q-1}_{\varsigma}, \infty\} 
\subset \CW_w
 \longrightarrow \ZZ_q, 
\ \ \
\left\{ 
\begin{array}{rcl} 
\varsigma & \longmapsto &  \varsigma
\\
\infty & \longmapsto & q .
\end{array}
\right.
$$

\begin{figure}[h!]
\begin{center}
\includegraphics[width=\textwidth]{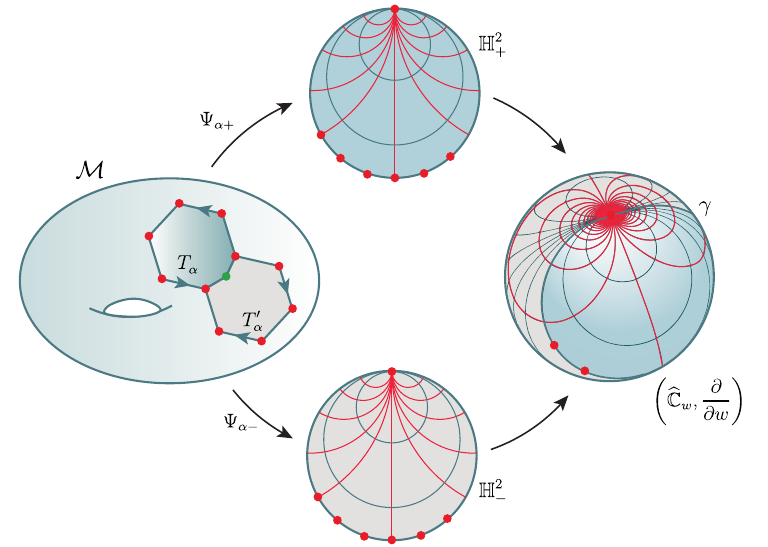}
\caption{
Two topological $q$--gons $T_1, \,T_{1}'$ are 
mapped to the
half planes $\HH_{+}, \, \HH_{-}$. 
We sketch the case  $q=7$,  
the critical points are red points in $\mathcal{M}$, 
a cocritical point is colored green. 
}
\label{fig:prueba-klein}
\end{center}
\end{figure}

\begin{remark}
\label{re:campo-esfera}
Throughout the work $(\CW_w, \partial / \partial w )$ 
denotes the Riemann sphere furnished with the 
holomorphic vector field $\partial / \partial w$, see Figure
\ref{fig:prueba-klein}. 
By abuse of notation,
this pair must be understood as the euclidean or
flat Riemannian metric, with a 
singularity at $\infty \in \CW_w$. 
The concepts of euclidean segments  and 
trajectories of the real vector field 

\centerline{$\Re{\e^{i \theta} \del{}{w}} \doteq 
\cos(\theta)\del{}{x}+
\sin(\theta)\del{}{y}$, \ \ \ $\theta \in [0, 2\pi]$,
}

\noindent
in $\CW_z$ are used in an equivalent way.
On Riemann surfaces, 
there exists a fruitfully relationship between complex analytic 
functions
and vector fields, it is described as the ``Dictionary'' 
\cite[prop.\,2.5]{AlvarezMucino3}.   
\end{remark}

\noindent 
Since each tile 
$T_\alpha$ or $T'_\alpha$ of $\mathscr{T}(\widehat{\Gamma})$
is an open Jordan domain,  there exist $C^1$ diffeomorphisms
\begin{equation}
\label{ec:cartas-en-las-losetas}
{\Psi}_{\alpha+}: 
T_\alpha \subset \mathcal{M} \longrightarrow 
\HH_+ \subset \left( \CW_w, \del{}{w} \right),
\ \ \
{\Psi}_{\alpha -}: 
T_{\alpha}^\prime \subset \mathcal{M} \longrightarrow 
\HH_- \subset \left( \CW_w, \del{}{w} \right),
\end{equation}

\noindent 
which extend continuously to the boundary 
$\partial \overline{T_\alpha}$, $\partial \overline{T'_\alpha}$. 
That is, 
$\Psi_{\alpha +}: \partial \overline{ T_\alpha } \longrightarrow
\gamma$ is a homeomorphism and the same property holds for 
$\Psi_{\alpha -}$.

\noindent 
The geometrical structure provided by the 
holomorphic vector field
$\del{}{w}$ on $\CW_w$ is as follows.
Recall that the real vector fields 

\centerline{$
\Re{\del{}{w}} \doteq \del{} {x}, 
\ \ \
\Im{\del{}{w}} \doteq  \del{}{y}
$}

\noindent
are real analytic on $\CW_w$ with a ``dipole'' at 
$w=\infty$
as the unique singularity.
The usual complex structure $J$ of $\CW_w$ is provided
by 

\centerline{$
J(\del{}{x}) = \del{}{y}, 
\ \ \ 
J(\del{}{y})= - \del{}{x}$.}

\noindent 
The pullback complex structure on 
$T_\alpha, J_{\alpha +}$ is

\centerline{$
J_{\alpha +} \big( (\Psi_{\alpha +})^*\del{}{x} \big) = 
(\Psi_{\alpha +})^* \del{}{y}, 
\ \ \ 
J_{\alpha +} \big( (\Psi_{\alpha +})^*\del{}{y} \big) = 
- (\Psi_{\alpha +})^* \del{}{x}
$.} 

\noindent 
Summing up, 
the pair $(T_\alpha, J_{\alpha +})$ is an open 
Riemann surface with 
holomorphic chart $\Psi_{\alpha +}$ as in
\eqref{ec:cartas-en-las-losetas}.
Analogously, 
each $(T_{\alpha}^\prime, J_{\alpha -})$ 
is an open Riemann surface.

{\it Step 3.} 
By using the labelled graphs
$(\widehat{\Gamma}, 
\mathcal{L}_c )$ and $(\gamma, \mathcal{L}_\gamma)$,
we construct a local homeomorphism

\centerline{
$\mathcal{R}: \mathcal{M} \backslash V (\widehat{\Gamma}) 
\longrightarrow \CW_w \backslash \{1, \ldots, q-1, \infty \} $,}

\noindent 
defined in the interior of the tiles as 
$\Psi_{\alpha +}$ or $\Psi_{\alpha -}$. 
Hence $\mathcal{R}$ satisfies that 
$\mathcal{R}(T_\alpha) =\HH_+$ and 
$\mathcal{R}(T_{\alpha}^\prime ) =\HH_-$.   

\noindent 
Let $\upzeta_0 \in E(\widehat{\Gamma}) \subset \mathcal{M}$ be
a point in an edge of $\widehat{\Gamma}$
but not a vertex of it, we want to construct a holomorphic chart
for $\upzeta_0$.  
There exist a unique trajectory of 
the $C^1$ vector field $(\Psi_{\alpha +}) ^* \del{}{y}$ in $T_\alpha$, 
say the vertical line $\{ x= x_0 \} \subset \HH_{+}$, 
such that $\upzeta_0$ is its unique extreme point. 
Figure \ref{fig:prueba-klein} illustrates this behavior.
We define  

\centerline{
$\mathcal{R}(\upzeta_0) = x_0 + 0i 
\in \RR \backslash \{1, \ldots, q-1, \infty \} \subset \CW_w$.}

\noindent 
Moreover, since $\upzeta_0$ is in the boundary of a gray tile, 
say $T_{\alpha}^\prime$, then 
there exists a unique trajectory of 
the $C^1$ vector field $(\Psi_{\alpha -})^* \del{}{y}$ in $T'_\alpha$, 
given by the vertical line $\{ x=x_0 \} \subset \HH_{+}$, with 
$\upzeta_0$ is its unique extreme point. 
Hence, $\mathcal{R}$ is well defined and continuous at 
$\upzeta_0$.  

{\it Step 4.} 
We construct a holomorphic chart that covers $\upzeta_0$.
Let $\epsilon$ the minimal euclidean distance of 
$\mathcal{R}(\upzeta_0)=x_0$ to the vertices of $\gamma$, thus
$\epsilon = 
\min \{ \vert x_0 - m  \ \vert \ m \in 1, \ldots, q-1   \} > 0$.
We denote by $D(x_0 , \epsilon ) \subset \CC$ 
the complex disk, 
with center $x_0 \in \CC$ and radius $\epsilon$. 
Define a continuous chart
\begin{equation}
\label{ec:carta-de-frontera}
\Psi_{\upzeta_0 }^{-1}: D(x_0 , \epsilon ) \subset \CC
\longrightarrow \mathcal{V} (\upzeta_0) \subset \mathcal{M}
\end{equation}

\noindent 
as follows:

\noindent 
i)
If $x + i0 \in D(x_0 , \epsilon ) \subset \CC \cap \RR$, then 
$\Psi_{\upzeta_0 }^{-1}(x)$ is a orientation
preserving isometry map from the euclidean segment
$(x_0 - \epsilon, x_0 + \epsilon)$
to the edge of $\widehat{\Gamma}$ that contains $\upzeta_0$.

\noindent 
ii)
If $x + iy \in D(x_0 , \epsilon )  \subset \CC \cap \HH_+ $, 
then 
$\Psi_{\upzeta_0 }^{-1}(x+y)\in T_\alpha$ is the unique point
such that under the real flow of the real vector field
$(\Psi_{\alpha +})^* \del{}{y}$ by time $-y$ arrives to 
$\Psi_{\upzeta_0 }^{-1}(x)$.

\noindent 
iii) 
If $x + iy \in D(x_0 , \epsilon ) \subset \CC \cap \HH_-$, 
then 
$\Psi_{\upzeta_0 }^{-1}(x+y)\in T_\alpha$ is the unique point
such that under the real flow of the real vector field
$(\Psi_{\alpha -})^* \del{}{y}$ by time $y$ arrives to 
$\Psi_{\upzeta_0 }^{-1}(x)$.

\begin{figure}[h!]
\centering
\includegraphics[width=0.45\textwidth]{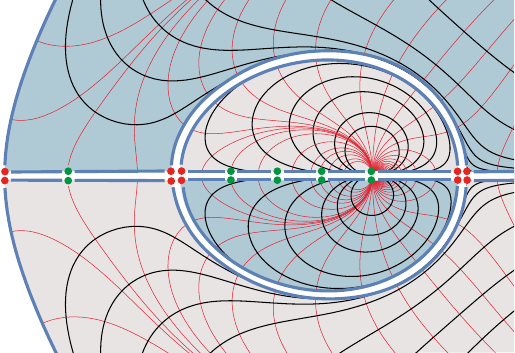}
\caption{
Zoom of Figure \ref{fig:tres-mosaicos-etiquetas-h}
illustrates four tiles before the gluing process.
The tiles are provided
with the real trajectories of the vector fields
$\Psi_{\alpha \pm}^* \del{}{x}$ in black and 
$\Psi_{\alpha \pm}^* \del{}{y}$ in red.
Moreover, the critical points are red, 
the cocritical points are green. 
A simple pole  of $R(z)$ appears as cocritical point.
}
\label{fig:prueba}
\end{figure}

\noindent 
A priori, the map 
$\Psi_{\upzeta_0}^{-1}$ is not $C^1$ at the
edges of $\widehat{\Gamma}$, when you consider the original 
$C^1$ structure on $\mathcal{M}$.

\noindent 
According to the 
Equations \eqref{ec:cartas-en-las-losetas} and
\eqref{ec:carta-de-frontera}, when the compositions 
$\Psi_{\alpha +}\circ \Psi_{\upzeta_0}^{-1}$, 
$\Psi_{\alpha -}\circ \Psi_{\upzeta_0}^{-1}$
are well defined, they are the restrictions of translations 
in $\CC$.
Hence, the collection 

\centerline{
$\{ \Psi_{\alpha +},\,   \Psi_{\alpha -}, \,  \Psi_{\upzeta_0} 
\ \vert \ 
\alpha \in 1, \ldots, n, \
\upzeta_0 \in E(\widehat{\Gamma}) \backslash V(\widehat{\Gamma})
\}$, }

\noindent  
is a holomorphic atlas for $\mathcal{M} \backslash V(\widehat{\Gamma})$. 
We denote the resulting open Riemann surface as $M^\circ$.

{\it Step 5.}
The construction of a homolomorphic map from 
$$
R^\circ: M^\circ \longrightarrow \CW_w
$$
is immediate. Consider the blue tile $T_1$
and the point $z_0 \in T_1 $ 
such that $\Psi_{1+}(z_0) = i$. 
We define $R^\circ = \Psi_{1+}$ in $T_1$.
If $T_{\beta}^\prime$ is any gray tile adjacent to $T_1$, 
then we extend our function $\Psi_{1+}$ 
by analytic continuation, that is 
$$
R^{\circ}: T_1  \cup T_{\beta}^\prime  \cup  \overline{z_j z_\iota}
\longrightarrow \CW_w, 
\ \ \
\mathfrak{z} \longmapsto 
\left\{
\begin{array}{lcl}
\Psi_{1+}(\mathfrak{z}) & \ \ \ \  &\mathfrak{z} \in T_1
\\
\Psi_{\beta -}(\mathfrak{z}) & &\mathfrak{z} \in T_{\beta}^\prime
,
\end{array}
\right. 
$$
where $\cup \, \overline{z_j z_\iota}$ 
denotes the interior of the edges
of $\widehat{\Gamma}$ in 
$\partial T_1 \cap T_{\beta}^\prime$.
In fact, with the above definition $\Psi_{1+}$
extends holomorphically to $M^\circ$, 
it produces a holomorphic function $R^\circ$, as we asserted.

{\it Step 6.}
Using the Riemann's extension theorem, we have that 
the vertices $V(\widehat{\Gamma})$ are
conformal punctures of $M^\circ$, 
since the above atlas can be
extended to $V(\widehat{\Gamma})$. 
Thus $M$ is a compact Riemann surface.
Let 

\centerline{$R: M \longrightarrow \CW_w$ }

\noindent 
be the analytic extension of 
$R^\circ$ to the vertices $V(\widehat{\Gamma})$.
We describe the local behavior of $R$
at each vertex. 
Note that, $R^\circ$ near $z_\iota$ 
is topologically equivalent to 

\centerline{
$\{ z \longmapsto z^{\nu/2}\}$, \ where  
$\nu=2$ or $\nu \geq 4$ and even,
}

\noindent 
means the valence of the vertex\footnote{
Recall that the vertices of $\widehat{\Gamma}$ split into two families: those of valence 2 and those of valence $\geq 4$, 
denoted as $\{ \zeta_\kappa \}$ and $\{ z_\iota \}$ respectively (see Step 1 and Step 3 of Schwarz--Klein's 
algorithm in Section \ref{sec:algoritmo-Schwarz-Klein-racional}). 
In what follows we will use the notation $\{ z_\iota \}$ for both families.
}
$z_\iota$, 
as a consequence:

\noindent
If $z_\iota \in V(\widehat{\Gamma})$ has label 
$j \in \{1, \ldots, q-1\}$ 
then $z_\iota$ is 

\quad\quad
$
\begin{cases}
\text{a critical point of } R \text{, with finite critical 
value, and ramification index } \nu/2 \text{ when } \nu \geq 4,\\
\text{a cocritical point of } R \text{ when } \nu=2.
\end{cases}
$

\noindent
If $z_\iota \in V(\widehat{\Gamma})$ has label 
$q$  
then $z_\iota$ is 

\quad\quad
$
\begin{cases}
\text{a pole of order } \nu/2 \text{ of } R \text{ when } \nu \geq 4,\\
\text{a simple pole of } R \text{ when } \nu=2.
\end{cases}
$

\noindent 
The construction of the rational function $R$ and 
assertion (2)  of Theorem \ref{teo:principal} 
are done.

\begin{corollary}
\label{cor-cubierta-ramificada}
\begin{enumerate}[label=\arabic*),leftmargin=*]
\item  
A ${\tt t}$--graph $\Gamma$ in $\mathcal{M}$, provided
with a consistent labelling $\mathcal{L}_c$,
determines a ramified cover 
$\mathcal{R}: \mathcal{M} \longrightarrow \CW_w$.

\item 
Conversely, a $C^1$ ramified cover  
$\mathcal{R}$  determines a (non unique)
rational function $R: M \longrightarrow \CC_w$.

\item
Let  $\mathscr{T}_\gamma(R)$ be a tessellation
of $M$, from
a rational function $R$ and a Jordan path $\gamma$. 
If the tiles  $T_\alpha, \, T_\alpha^\prime$ are
adjacent along an edge  $\overline{z_j z_\iota}$
of $\widehat{\Gamma}$, 
then there is an branch of the inverse function

\centerline{$ 
R^{-1}: 
\big( \CW_w \backslash \gamma \big)
\cup R(\overline{z_j z_\iota} )
\longrightarrow 
T_\alpha \cup \overline{z_j z_\iota} \cup T_\alpha^\prime
\subset M
$,}

\noindent 
where $\big( \CW_w \backslash \gamma \big)
\cup R(\overline{z_j z_\iota} )$ is a maximal domain of univalence.   
\hfill
$\qed$
\end{enumerate}
\end{corollary}

Some comments are in order. 

The function $R$ in Theorem \ref{teo:principal} 
assertion (2) is not 
unique, due to several choices:

\noindent 
i) The complex structure on $M$ 
depends on the choice of the critical values 
$\{ 1, \ldots, q-1, \infty\} \subset \CW_w$

\noindent 
ii) For any  $q$ distinct critical values 
$\{ w_1, \ldots , w _q\} \subset \CC_w$, the construction with the same technique remains simple
when $\gamma$ is a polygonal path in $\CW_w$, 
as we will show in 
Corollary \ref{cor:funciones-racionales-por-cirujia}.

\noindent 
iii) 
For
fixed $\{ w_1, \ldots , w _q\}$ the change of $\gamma$
was studied by Habsch \cite{Habsch}, see also 
\cite[\S\,5.4]{Lando-Zvonkin}; both in the case $M=\CW_z$.
However, note that if $\gamma$ is allowed to vary by an isotopy relative to the $q$ 
distinct critical values, then the function $R$ is unchanged.

\subsection{Relationship between analytical and combinatorial structures}

As a matter or record, let us consider the following.

\noindent 
$\bigcdot$
A homogeneous tessellation $\mathscr{T}_\gamma (R)$, 
or equivalently an  analytical triple  $(M, R, \gamma)$, denotes a rational function
$R: M \longrightarrow \CW_w$, 
on a compact Riemann surface, and 
a cyclic order $\mathcal{L}_\gamma$ of its
$q$ critical values $\mathcal{CV}_{R}$ of $R$. 

\noindent
$\bigcdot$ 
A pair $( \widehat{\Gamma}, \mathcal{L}_c)$, where 
$\widehat{\Gamma} \subset \mathcal{M}$ is an ${\tt R}$--map
(on a $C^1$ compact oriented surface $\mathcal{M}$),
with vertices of even valence $\geq 2$, 
and a consistent $q$--labelling  $\mathcal{L}_c: V(\widehat{\Gamma})
\longrightarrow \ZZ_q$.

\noindent 
$\bigcdot$
A pair $(\Gamma, \mathcal{L}_c)$, 
where $\Gamma \subset \mathcal{M}$ is a ${\tt t}$--graph,
with vertices of even valence $\geq 4$, and 
a consistent $q$--labelling  $\mathcal{L}_c: V(\Gamma)
\longrightarrow \ZZ_q$.

\noindent
$\bigcdot$ 
$\mathscr{T}=\mathscr{T}(\Gamma) = \mathcal{M} \backslash \Gamma $ is a 
tessellation.

\noindent 
Note that only $\mathscr{T}_\gamma(R)$
belongs to the complex analytic category. 
The 
$q$--labellings of $\widehat{\Gamma}$ and $\Gamma$
coincide on the vertices of $\Gamma$. 
We have the following correspondences:

\begin{center}
\begin{picture}(160,78)

\put(-180,4){\vbox{\begin{equation}
\label{dia:correspondencia-completa}\end{equation}}}

\put(-86,4){$\mathscr{T}(\Gamma)$}

\put(52,65){$(\widehat{\Gamma}, \mathcal{L}_c)$}

\put(-50,7){\vector(1,0){38}}
\put(-20,7){\vector(-1,0){38}}
\put(-38,13){$\textcircled{{\bf{\tiny 1}}}$}
\put(-7,4){$\Gamma$}

\put(105,74){$\textcircled{{\bf{\tiny 5}}}$}
\put(114,25){$\textcircled{{\bf{\tiny 6}}}$}

\put(48,7){\vector(-1,0){45}}
\put(22,13){$\textcircled{{\bf{\tiny 2}}}$}
\put(52,4){$(\Gamma, \mathcal{L}_c)$}

\put(64,59){\vector(0,-1){45}}
\put(51,34){$\textcircled{{\bf{\tiny 3}}}$}

\put(68,14){\vector(0,1){45}}
\put(70,34){$\textcircled{{\bf{\tiny 4}}}$}

\put(134,65){$ \mathscr{T}_\gamma(R).$}
\put(87,7){\vector(1,1) {52}}
\put(130,67){\vector(-1,0){45}}

\end{picture}
\end{center}

\noindent
Map
$\textcircled{{\bf{\tiny 1}}}$ is an equivalence,
by definition.

\noindent 
Map $\textcircled{{\bf{\tiny 2}}}$ forgets the 
labels
$\{\mathcal{L}_c(z_\iota) \}$ of
the vertices of $\Gamma$. 

\noindent 
Map $\textcircled{{\bf{\tiny 3}}}$ forgets the vertices 
$\{\upzeta_\kappa \} \subset \widehat{\Gamma}$ of valence $2$
and their labels
$\{\mathcal{L}_c(\upzeta_\kappa) \}$. 

\noindent 
Map $\textcircled{{\bf{\tiny 4}}}$ is the edge subdivision 
of $\Gamma$, that introduces vertices of valence 2, according 
to the consistent $q$--labelling $\mathcal{L}_c$ of $V(\Gamma)$.

\noindent 
Map $\textcircled{{\bf{\tiny 5}}}$ 
is assertion (1) in Theorem \ref{teo:principal}, 
where $\widehat{\Gamma} = R^* \gamma$
and the consistent $q$--labelling
$\mathcal{L}_c : V(\widehat{\Gamma}) \longrightarrow \ZZ_q$ 
is the pullback $R^*\mathcal{L}_\gamma$ 
of the cyclic order $\mathcal{L}_\gamma$
of the critical values of $R$.

\noindent 
Map $\textcircled{{\bf{\tiny 6}}}$ 
is assertion (2) in Theorem \ref{teo:principal}.

\noindent 
The tessellations 
$\mathscr{T}(\widehat{\Gamma})$ 
and $\mathscr{T}(\Gamma)$ coincide,
set theoretically.

\smallskip

The construction of $R$ in the proof of 
Theorem \ref{teo:principal} produces critical
values in $\RR \cup \{ \infty \}$. 
We introduce now a more flexible construction.

\begin{corollary}[Construction of rational functions by surgery]
\label{cor:funciones-racionales-por-cirujia}
Let 
$\overline{\mathcal{P}} \subset (\CW_w, \partial / \partial w )$ 
be a closed polygon with vertices
$\{ w_1, \ldots , w_q\} \subset \CW_w$ and 
oriented euclidean segments as edges, 
where
$\mathcal{P}$ is on the left side of the edges. 
Assume that $\mathcal{P}$ and its 
complement $\overline{\mathcal{P}}^c$ are open Jordan domains.
Consider a collection of $2n$ tiles
$$
T_\alpha = \Big(\mathcal{P}, \del{}{w} \Big) ,
\ \ \  
T_\alpha^\prime = 
\Big(\overline{\mathcal{P}}^c, \del{}{w} \Big),
\ \ \  
\alpha \in \{1, \ldots ,n \geq 2\}.
$$

\noindent 
Assume that $M$ is a compact Riemann surface obtained 
by gluing together closed tiles $\overline{T_\alpha}$ 
and $\overline{T_\beta^\prime}$, along an edge
by euclidean isometries. 
Each vertex $z_\iota \in M$ (arising from the gluing of  vertices 
in the closed tiles)
has cone angle $2\pi \nu_\iota$, $\nu_\iota \in \NN$. 
Then $M$ determines
\begin{enumerate}[label=\roman*),leftmargin=*]
\item 
a rational function $R: M \longrightarrow \CW_w$ of degree $n$
with tessellation 
$\mathscr{T}_{\gamma} (R)$, 
where $\gamma =\partial \overline{\mathcal{P}}$, and 

\item
a pullback rational vector field $X = R^* \del{}{w}$ on $M$. 
\end{enumerate}
\end{corollary}

The case $n=\infty$ gives origin to complex analytic functions
(not necessarily meromorphic) on non compact Riemann surfaces, 
see \cite{AlvarezMucino5}.

\begin{proof} 
The boundaries of the tiles determine an ${\tt R}$--map on $M$. 
Note that, the alternating color condition for the tiles
of $\mathscr{T}_{\gamma} (R)$ is implicit in the hypothesis.
The set of points  $\{ w_1, \ldots , w_q\}$ contains 
the critical value set $\mathcal{CV}_R \neq \emptyset$. 
If no minimality condition is assumed on the edges of 
$\overline{\mathcal{P}}$, 
then additional subdivision of an edge produces fake
cocritical points and fake critical values. 
Furthermore,  a straight line segment 
through $\infty$ is allowed
as an edge in $\partial \overline{\mathcal{P}}$. 
The pullback  of $\partial / \partial w$
under $R(z)$ is 
$(1/R^\prime(z)) \frac{\partial }{\partial z}$, 
this analytic expression makes sense 
on any Riemann surface. 
\end{proof}

\begin{corollary}
\label{cor:el-orden-ciclico-para-los-valores}
Let $R:M \longrightarrow \CW_w $ be a rational function.
Topologically
the tessellation $\mathscr{T}_\gamma (R)$ 
only depends  
on the cyclic order  $\mathcal{L}_\gamma$,
given to the critical values
$\{ w_1, \ldots , w_q \}$ of $R$.
\hfill\qed
\end{corollary}

In other words, 
the tessellation $\mathscr{T}_\gamma (R)$ 
is topologically independent  
on the particular shape of $\gamma$
(as long as $\gamma$ is in the same isotopy class relative to 
the $q$ distinct critical points). 
There are $q!$ cyclic orders $\mathcal{L}_\gamma$ for 
the critical values $\{ w_1, \ldots, w_q \}$.
Note that for  $q\geq 4$, up to M\"obius transformation,
the usual choice for the critical values is $\{0,1, \infty, w_4, \ldots, w_q \}$.
A complete discussion is to appear in \cite{AlvarezMucino5}.

\section{Examples}
\label{sec:Examples}
Given integers ${\tt p}$, ${\tt q}$, ${\tt r} \geq 2$,
let $G_{{\tt p},{\tt q},{\tt r}}$ be the group with generators 
and relations

\centerline{$
\left\{
g_1, \, g_2, \, g_3 
\ \vert \
g_1^2 = g_2^2 =g_3^2=1,
\ \
(g_1 g_2)^{\tt p} =(g_2 g_3)^{\tt q} =(g_3 g_1)^{\tt r}=1
\right\} $ . }

\noindent 
Then $G_{{\tt p},{\tt q},{\tt r}}$ can be represented 
as a group of isometries of one of the 2--dimensional geometries:
\begin{equation}
\label{ec:tres-geometrias}
\dfrac{1}{{\tt p}} +\dfrac{1}{{\tt q}} +\dfrac{1}{{\tt r}}
\
\left\{
\begin{array}{lcl} 
>  1  & \ \  & \mbox{ spherical } \mathbb{S}^2 ,
\\ 
= 1  && \mbox{ euclidean } \RR^2 ,
\\
< 1  && \mbox{ hyperbolic } \HH \, .
\end{array}
\right.
\end{equation}

\noindent  
There exists a triangle 
$T= T_{{\tt p},{\tt q},{\tt r}}$ in the respective space
$\Omega_z$ ($\mathbb{S}^2$, $\RR^2$ or $\HH$) bounded by geodesics edges and 
with interior angles ${\pi}/{\tt p}$, 
${\pi}/{\tt q}$, ${\pi}/{\tt r}$.

\begin{theorem*}
The group $G_{{\tt p},{\tt q},{\tt r}}$  is isomorphic to 
the group of isometries of $\Omega_z$,
generated by reflections on the 3 geodesic edges of 
$T \subset \Omega_z$. 
The triangle $T$ is a fundamental domain for the action of 
$G_{{\tt p},{\tt q},{\tt r}}$  on $\Omega_z$.
\hfill
\qed
\end{theorem*}

Let $H_{{\tt p},{\tt q},{\tt r}}$ 
be the index 2 subgroup of $G_{{\tt p},{\tt q},{\tt r}}$,
its elements are holomorphic automorphisms of 
$\Omega_z$.

\begin{corollary}
\begin{enumerate}[label=\arabic*),leftmargin=*]
\item
The action gives rise to a tessellation 

\centerline{
$\mathscr{T}_{{\tt p},{\tt q},{\tt r}} =
\{g(T) \ \vert \ g \in G_{{\tt p},{\tt q},{\tt r}}\}
$
}

\noindent 
of $\Omega_z$, with blue tiles 
$\{g(T) \ \vert \ g \in H_{{\tt p},{\tt q},{\tt r}}\}$.

\item
There exist meromorphic functions
$R:\Omega_z \longrightarrow \CW_w$,
each of degree the order of $H_{{\tt p},{\tt q},{\tt r}}$ 
and $q=3$ critical values.
Each function $R$ and its tessellation
$\mathscr{T}_{{\tt p},{\tt q},{\tt r}}$
are $H_{{\tt p},{\tt q},{\tt r}}$--invariant. 
\hfill
\qed
\end{enumerate}
\end{corollary}

\begin{example}[Schwarz's tessellations]
As far as we known, 
the first examples of rational functions that origin tessellations 
are due to H.\,A.\, Schwarz \cite{Schwarz}.
For each ${\tt p}, \, {\tt q}, \, {\tt r}$,
as in Equation \eqref{ec:tres-geometrias},
there exists tessellation $\mathscr{T}_{{\tt p},{\tt q}, {\tt r}}$ by 
triangles (tiles) $\{ T_\alpha, \, T_\alpha^\prime\}$ in $\mathbb{S}^2$, see Figure \ref{fig:solidos-platonicos}.

\begin{figure}[h!]
\begin{center}
\includegraphics[width=0.8\textwidth]{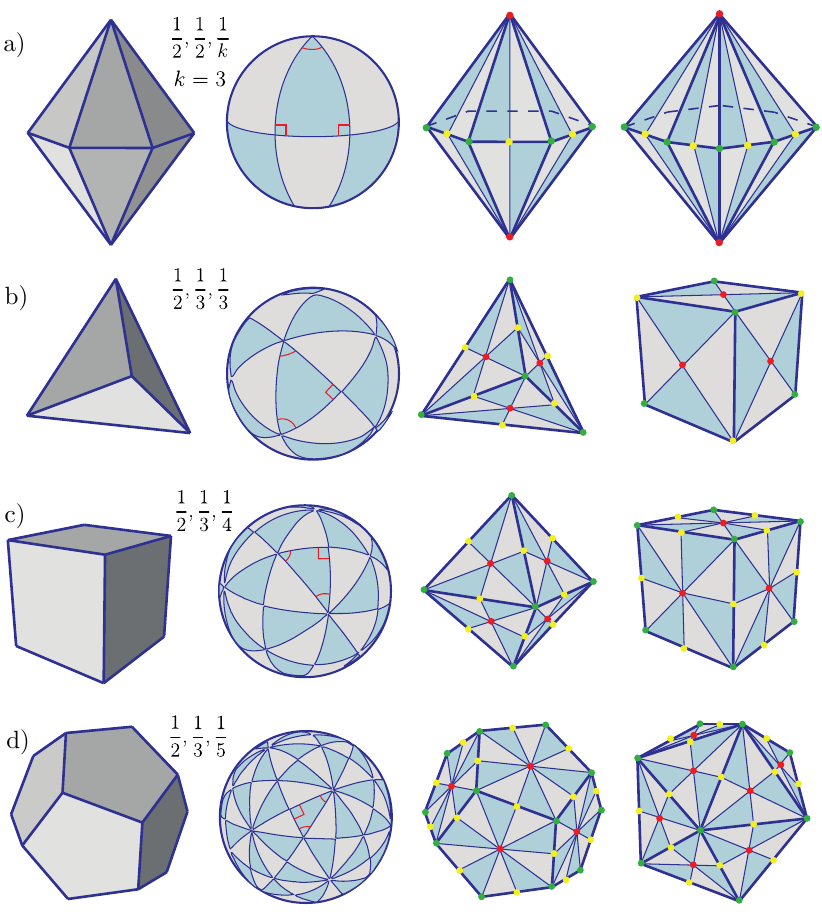}
\caption{
Schwarz's tessellations.
The left column illustrates the $k$--dihedron and the platonic solids. 
The other columns show tessellations of rational 
$G$--invariant Bely{\u \i} functions.
}
\label{fig:solidos-platonicos}
\end{center}
\end{figure}

\noindent 
Let $\mathfrak{P}$ be a $k$--dihedron or a platonic solid, 
see Figure \ref{fig:solidos-platonicos} first column.
We make the following conventions.
The vertices  
$\{ \mathfrak{v}_\iota\} \in \mathfrak{P}$
(green points) are zeros of $R$;
the face centers  ${\tt r} \in \mathfrak{P}$
(red points) are poles of $R$;
the middle points  
$\{ \mathfrak{e}_\kappa\}$ of edges
of $\mathfrak{P}$ (yellow points) 
assume the value $1$. 

\noindent
Moreover,
we assign suitable multiplicities $\mu, \nu \geq 2$ to the 
vertices $\{v_\iota\}$ and the 
centers $\{ c_j\}$,
according to ${\tt p},\, {\tt q}, \, {\tt r}$. 
The $H_{{\tt p} {\tt q} {\tt r}}$--invariant functions are  

\centerline{$
R(z)= 
\lambda
\dfrac{\prod\limits_\iota (z-\mathfrak{v}_\iota)^{\mu_\iota} }{
\prod\limits_j (z-\mathfrak{c}_j)^{\mu_j}  },
\ \ \ \ \
\left\{
\begin{array}{l}
R(\mathfrak{v}_\iota)=0 ,
\\
R(\mathfrak{c}_j) = \infty ,
\\
R(\mathfrak{e}_\kappa) = 1\, .
\end{array}
\right.
$
}

\noindent 
They are Bely{\u \i} functions with critical values $\{0, 1, \infty \}$.
The explicit formulas for the $G$--invariant functions
are cumbersome, for the icosahedron, we get 
$P(z)/Q(z)$ with

\centerline{$
P(z)=z^5 \left(z^5-\left(1-\sqrt{5}\right)^5\right)^5 \left(z^5-\left(1+\sqrt{5}\right)^5\right)^5 
\ \ \mbox{and } 
$}

\centerline{$
\begin{array}{l}
Q(z)=\left(z^5-\left(\frac{1}{4}
\left(-1+\sqrt{5}+\sqrt{6 \left(5+\sqrt{5}\right)}\right)+\frac{2 i}{\sqrt{3} \left(1-\frac{1}{\sqrt{15+6 \sqrt{5}}}\right)}\right)^5\right)^3 
\\
\vspace{-.3cm}
\\
\ \ \ \ \ \ \ \, \left(z^5-\left(-\frac{2 \sqrt{3} \left(2+\sqrt{5}\right)}{\sqrt{3}+3 \sqrt{5+2 \sqrt{5}}}-i\frac{2}{\sqrt{3} \left(1+\frac{1}{\sqrt{15+6 \sqrt{5}}}\right)}\right)^5\right)^3 
\\
\vspace{-.3cm}
\\
\ \ \ \, \left(z^5-\left(-\frac{\sqrt{30 \left(5+\sqrt{5}\right)}}{15+\sqrt{75+30 \sqrt{5}}}-i\frac{5 \sqrt{18-6 \sqrt{5}}}{15+\sqrt{75+30 \sqrt{5}}}\right)^5\right)^3
\\
\vspace{-.3cm}
\\
\left(z^5-\left(\frac{1}{2} \left(2+\sqrt{5}+\sqrt{15+6 \sqrt{5}}\right)+i \sqrt{\frac{1}{2} \left(10+\sqrt{5}+\sqrt{75+30 \sqrt{5}}\right)}\right)^5\right)^3.
\end{array}
$
}
\end{example}

\begin{example}[Fortunate rational functions]
We say that a rational function $R(z)$ is {\it fortunate}
when its critical values belong to 
$\RR \cup \{ \infty \}$, in this case a natural choice for $\gamma$
is this circle itself. 
Hence, a fortunate rational function $R$ has natural tessellation 

\centerline{
$M \backslash \Gamma$,
\ where $\Gamma = \{ \Im{R(z)} =0 \}$. 
}

\noindent 
Some advantages of the tessellations of fortunate rational functions are: 

\noindent 
$\bigcdot$
The graph $\Gamma$ is parametrized by  
the solution of the equations
$\{ R(z)=c \  \vert \ c \in \RR \}$.

\noindent 
$\bigcdot$
For $M=\CW_z$, as fas as we known $\Gamma$ is 
a reducible real algebraic curve
of degree $n$ (the degree of $R(z)$). 

\noindent
Let $\ZZ_p$ be the cyclic $p$--group, acting on $\CC_z$ by rotations  
$\{ r_s(z) \ \vert \ s=1, \ldots, p\}$.
It is possible to construct fortunate functions with
$\ZZ_p$--symmetry, as follows. 
For $p \geq 2$,
we use the family of rational functions

\centerline{$R(z)= \int_0 ^z  (P(\zeta)/Q(\zeta) ) d\zeta$} 

\noindent 
having the following properties:

\noindent 
$\bigcdot$
The critical points $\{ z_\iota \}$ of $R(z)$
with their multiplicities are $\ZZ_p$--invariant.

\noindent
$\bigcdot$ 
The coefficients of $R(z)$ are real.

\noindent
$\bigcdot$ 
For each critical point $z_\iota$, 
there exists a rotation $r_s(z)$ such that 
$r_s(z_\iota) \in \RR$ and its critical value
can be expressed as a real integral, thus \
$\int_0^{z_\iota} (P(\zeta)/Q(\zeta) )d\zeta
=
\int_0^{r_s(z_\iota)} (P(\zeta)/Q(\zeta) )d\zeta.
$

\noindent 
Hence, such $R(z)$ is fortunate.  
In the other direction, 
a $\ZZ_p$--invariant rational function is not
necessarily fortunate.
Figure \ref{fig:racional-inverso} illustrates the tessellations of the following fortunate functions
$$
\begin{array}{lr}
a)\ P(z) = 
\int_{0}^{z} \left( \zeta^5 (\zeta^3-1)^2 
(\zeta^3+1)^4 (\zeta^6-2^6)^3\right)  d\zeta, 
& \ \ \ 
q=5;
\\
& \vspace{-.3cm}
\\
b) \ R_1(z)=
\dfrac{z^5 (z^5-4^5)^5(z^5-5^5)^5}{(z^5-1^5)^3(z^5-2^5)^2 (z^5-3^5)}, 
& \ \ \ 
q=7;
\\ 
& \vspace{-.2cm}
\\
c)\ R_2(z)=
\dfrac{(z^5-2^5)^5(z^5-4^5)^5 }{z^5 (z^5+1^5)^2(z^5+3^5) (z^5+5^5)^3},
& \ \ \ 
q=7.
\end{array}
$$
\begin{figure}[h!]
\begin{center}
\includegraphics[width=\textwidth]{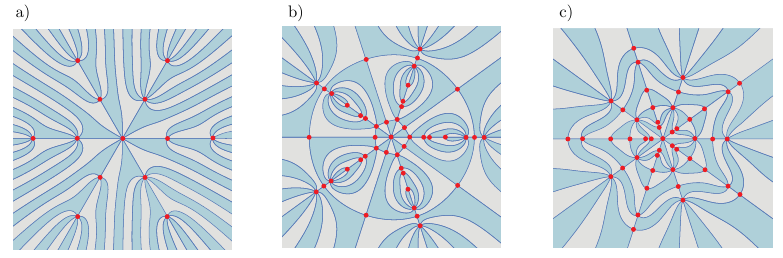}
\caption{
Affine tessellations of three fortunate functions. 
a) A  $\ZZ_3$--invariant polynomial whose tiles are $5$--gons.
(b)--(c)
Two $\ZZ_5$--invariant rational functions whose tiles are $7$--gons
(the cocritical points are not illustrated).
}
\label{fig:racional-inverso}
\end{center}
\end{figure}
\end{example}

\begin{example}
\label{ex:funcion-wp}
{\bf 1.}
 Weierstrass $\wp$--functions.
Let $\Lambda = \{ n + m \tau \ \vert \ 
(n,m)\in \ZZ^2, \Im{\tau } >0
\}$ be a lattice in $\CC_z$. 
The associated  Weierstrass $\wp$--function on the 
$\Lambda$--torus $M$ is 

\centerline{
$\wp(z): M = \CC_z/ \Lambda \longrightarrow \CW_w$. }

\noindent 
It has critical points  
$\mathcal{CP}_\wp = 
\{[0], \, [1/2], [\tau /2], \, [(1 + \tau )/2] 
\}$, 
with critical values 
$\mathcal{CV}_\wp= \{\infty, e_1= \wp(1/2), e_2, e_3 \}$, where 
$\{ e_k \}$ are the roots of the cubic equation
$\{ 4 e^3 - 60G_2 e - 140 G_3 =0 \}$,
for $G_k = \sum_{\omega  } (1/ \omega^{2k} )$, $\omega \in 
\Lambda \backslash \{ 0 \}$.
The ${\tt R}$--map $\widehat{\Gamma}$ has $q=4$,
Figure \ref{fig:cuadrilateros-2}.a describes the tessellations
on a fundamental domain of the torus, and 
the consistent $4$--labelling.    
The degree of $\wp(z)$ is $2$. 

\noindent 
Clearly, 
the ``chessboard'' in Figure \ref{fig:cuadrilateros-2}.a,
can be generalized to a chessboard with $2\sigma_1$ rows
and $2\sigma_2$ columns,
for a pair $(\sigma_1, \sigma_2) \in \NN^2$. Each pair
gives rise to a new torus 
$\widehat{M}= \CC_z / \Lambda_{\sigma_1,\sigma_2} $,
see Figure \ref{fig:cuadrilateros-2}.b.
The corresponding rational function $\widehat{R}$ on the torus
$\widehat{M}$ has degree $2\sigma_1 \sigma_2$.

\noindent
{\bf 2.} 
The tessellations of
$\CC_z$ by equilateral triangles, 
${\tt p}={\tt q}={\tt r} = 3$ in  Equation
\eqref{ec:tres-geometrias}, 
or more generally by two
non isometric triangles $T, \, T^\prime$,
are obtained from the derivatives
$\wp^\prime(z)$, see Figure \ref{fig:cuadrilateros-2}.c
They originate from elliptic functions 
of degree $3$ with $q=3$ critical values
on the corresponding 
torus $\CC_z/ \Lambda$.  

\noindent 
{\bf 3.} 
Functions on hyperelliptic surfaces.
Let $\widehat{\Gamma} \subset \mathcal{M}$ 
be
and ${\tt R}$--map 
on a surface $\mathcal{M}$ of genus $g \geq 1$, the
consistent $6$--labelling is shown in Figure 
\ref{fig:cuadrilateros-2}.d for genus 2. 
The tessellation $\mathcal{M} \backslash 
\widehat{\Gamma}$ is by four $(2{\tt g}+2)$--gons with 
$2{\tt g}+2$ vertices of valence 4. 
The associated 
rational function $R$ of degree 2 determines
an hyperelliptic Riemann surface  $M=(\mathcal{M}, J)$. 
There are
$2{\tt g} +2$ critical points of multiplicity 2 of $R$.
For each genus  ${\tt g}\geq 1$, these rational functions $R$
attain the minimum number of critical points in $\mathcal{M}$.
For an hyperelliptic Riemann surface

\centerline{
$M =\{ \mathfrak{z}^2 - (w-w_1) \cdots (w-w_{2{\tt g}+2}) =0 \}
\subset \CC^2=\{(w, \mathfrak{z})\}$, 
}

\noindent 
with 
$\{ w_j \ \vert \ j \in  1,  \ldots, 2{\tt g}+2 \}$
distinct points, 
the projection $\pi_1: M \longrightarrow \CC_w$
belongs to the above family of 
rational functions.

\begin{figure}[h!]
\begin{center}
\includegraphics[width=\textwidth]{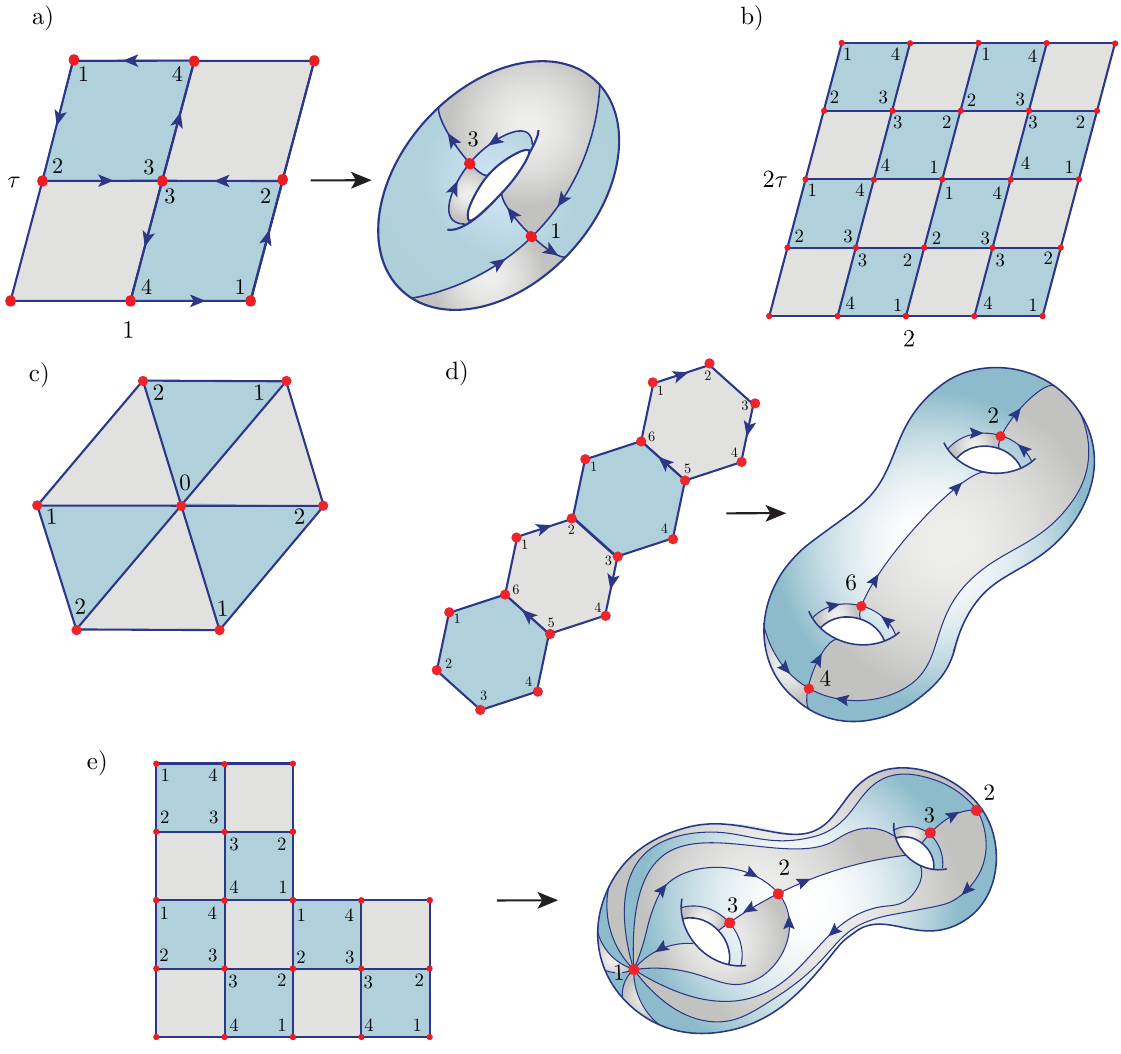}
\caption{
Tessellations on compact Riemann surfaces and their
consistent $q$--labellings.
a) A Weierstrass function $\wp(z)$ on a torus. 
b) A lift of $\wp(z)$ to a covering torus.
c) Derivative $\wp^\prime(z)$ on a torus.
d) A rational function on a hyperelliptic surface. 
e) A rational function on a surface of genus 2.
}
\label{fig:cuadrilateros-2}
\end{center}
\end{figure}

\noindent
{\bf 4.} 
Consider a ``L--chessboard'' 
as in Figure \ref{fig:cuadrilateros-2}.e
Under identification of the horizontal and vertical edges of the boundary, 
we get a surface $M$ of genus 2 
with a rational function $R: M \longrightarrow \CW_w$ of degree 6.
There are 9 critical points of multiplicity 2 and
1 critical point of multiplicity 6. 
This example generalizes to any genus ${\tt g} \geq 3$.
\end{example}

\subsection{Closing comments}
\label{sec:closing-comments}

\noindent
$\bigcdot$
In the case $\mathcal{M}= \mathbb{S}^2$, 
the existence of a consistent $q$--labelling for  
tessellations that arise from generic rational functions,  
was proved in  S.\,Koch {\it et al.} \cite{Koch-Lei};
by proving that the tessellations with
 $2n$ tiles and $2n-2$ vertices
of valence 4
satisfy global and local balance conditions. 
The existence of a consistent $q$--labelling for tessellations 
given by generic polynomials is in 
L.\,J.\,Gonz\'alez--Cely {\it et al.} \cite{Gonzalez-Mucino}. 
In the non generic case, the existence of the consistent $q$--labelling
remains as an open problem, even for $\mathcal{M}= \mathbb{S}^2$.

\noindent
$\bigcdot$ 
Given a tessellation $\mathscr{T}$ of $\mathcal{M}$
as in assertion (2) of Theorem \ref{teo:principal};
how can we provide an analytic description 
(in terms of coefficients for instance) of a corresponding
$R$ on $M$?

\noindent
$\bigcdot$
Consider a non homogeneous tessellation $\mathscr{T}$ on $\mathcal{M}$
that admits a consistent $q$--labelling $\mathcal{L}_{c}$, 
thus $\mathscr{T}$ is realizable from a rational function.
Several natural questions arise:

\noindent
How many different consistent $q$--labellings are there?

\noindent 
Under what conditions the value of $q$ is unique?

\noindent 
Example \ref{ex:racional-cocriticos} already shows
these phenomena.

\noindent
$\bigcdot$ 
Consider the space of analytic objects 
$\{ (M, R, \mathcal{L}_\gamma) \}$, 
where the genus of $M$, 
the degree of $R$ and the number of 
critical values $q$ are fixed. 
Considering the configurations 
of critical values $\{ w_1, \ldots , w_q \}$
of $R$ as free parameters, and recalling
that $\mathcal{L}_\gamma$ is the cyclic 
order in the critical values of $R$,
it is clear that $\mathcal{L}_\gamma$
is a discrete quantity. 
Very roughly speaking, 
using Diagram \ref{dia:correspondencia-completa}, 
we have the map  
$$
\textcircled{{\bf{\tiny 3}}} \circ 
\textcircled{{\bf{\tiny 5}}}: \{ (M, R,  \mathcal{L}_\gamma ) \} 
\longrightarrow \{(\mathcal{M}, \Gamma, \mathcal{L}_c, 
\{ w_1, \ldots , w_q \}) \}. 
$$

\noindent 
Thus we have a certain kind of fibration, with different topological fibers; from a family of analytic objects to topological data. 
The map $\textcircled{{\bf{\tiny 6}}}$ is a section of this fibration.  
An accurate formulation of these maps and their study
is a future project.

\noindent
$\bigcdot$ 
A meromorphic function $w:\CC_z \longrightarrow \CW_w$ is 
Speiser when it has a finite number of critical and asymptotic 
values, {\it i.e.} it belongs to the family $F(w_1, \ldots , w_q)$. 
The study of Speiser functions $w(z)$ on simply connected
Riemann surfaces was considered by 
A.\,Speiser \cite{Speiser} and 
R.\,Nevanlinna \cite[Ch.\,XI]{Nevanlinna2}. 
The study of tessellations and applications to
differential equations is the subject of 
A. Alvarez--Parrilla {\it et al.} \cite{AlvarezMucino5}.

\medskip

{\bf Acknowledgements.}

\noindent
The authors would like to thank 
Enrique Casta\~neda Alvarado, 
Fernando Orozco Zitli, 
Leidy J. Gonz\'alez--Cely;
for suggestions which have helped to improve the work.

%

\bibliographystyle{spmpsci}      


\end{document}